\documentclass[10pt]{article}
\usepackage[left=3cm,right=3cm,top=3cm]{geometry}

\usepackage[english]{babel}
\usepackage[utf8]{inputenc}
\usepackage{amsmath}
\usepackage{amssymb}
\usepackage{amsthm}
\usepackage[all,cmtip]{xy}
\usepackage[bookmarks=true]{hyperref}
\usepackage{tikz}
\usepackage{tikz-cd}
\usepackage{todonotes}
\usepackage{slashed}
\usepackage{mathtools}
\usepackage{dsfont}
\usepackage{tensor}
\usepackage{verbatim}
\usepackage{graphicx}
\usepackage[shortlabels]{enumitem}
\usepackage{xstring}\usepackage[capitalize]{cleveref}
\usepackage{tensor}
\usepackage{mathrsfs}
\usepackage{nicefrac}
 
\usepackage[mathscr]{euscript}

\definecolor{darkolivegreen}{rgb}{0.33, 0.42, 0.18} 
\definecolor{cobalt}{rgb}{0.0, 0.24, 0.43}

\newcounter{comments}
\newenvironment{displaycomment}{\begin{list}{}{\rightmargin=1cm\leftmargin=1cm}\item\sf\begin{small}}{\end{small}\end{list}}
\def\comments{
  \setcounter{comments}{1}
  
  }

\comments

\newcommand{\C}{\mathbb{C}}
\newcommand{\R}{\mathbb{R}}

\newcommand{\mc}[1]{\mathcal{#1}}
 
\newcommand{\Del}{\mathbb{D}}

\newcommand*{\defeq}{\mathrel{\vcenter{\baselineskip0.5ex \lineskiplimit0pt
                        \hbox{\scriptsize.}\hbox{\scriptsize.}}}%
        =}

\newcommand{\f}{\mathrm{f}}
\newcommand{\g}{\mathrm{g}}
\newcommand{\GL}{\operatorname{GL}}

\newcommand{\Hdiffhat}{\widehat{\mathrm{Conf}}}

\newcommand{\disk}{\mathbb{D}}

\newcommand{\sphere}{\overline{\mathbb{C}}}

\def\brackets#1{\IfStrEq{#1}{-}{}{(#1)}}

\theoremstyle{definition}
\newtheorem{definition}{Definition}[section]
\newtheorem{remark}[definition]{Remark}

\theoremstyle{plain}
\newtheorem{theorem}[definition]{Theorem}
\newtheorem{proposition}[definition]{Proposition}
\newtheorem{lemma}[definition]{Lemma}
\newtheorem{corollary}[definition]{Corollary}

\newtheorem{problem}[definition]{Problem}

\setlist{topsep=-0.2em,itemsep=0em}

\crefname{enumi}{\unskip}{\unskip}

\definecolor{darkolivegreen}{rgb}{0.33, 0.42, 0.18} 
\definecolor{cobalt}{rgb}{0.0, 0.24, 0.43}

\makeatletter
\hypersetup{
colorlinks=true,
urlcolor=cobalt,
linkcolor=cobalt,
citecolor=darkolivegreen,
pdfmenubar=false,
pdftoolbar=true,
pdftitle = {\@title},
}
\makeatother

\title{A Fermionic Grunsky operator}

\author{Peter Kristel, Eric Schippers \& Wolfgang Staubach}

\begin{document}

\maketitle

\begin{abstract}  To a conformal map $\f$ from  the disk $\disk$ into the complex plane onto a domain with rectifiable Ahlfors-regular boundary, we associate a new kind of Grunsky operator on the Hardy space of the unit disk. This is analogous to the classical Grunsky operator, which itself can be viewed as an operator on Bergman or Dirichlet space.  
  We show that the pull-back of the Smirnov space of the complement of $\f(\disk)$ by $\f$ is the graph of the Grunsky operator.  We also characterize those domains with rectifiable Ahlfors-regular boundaries such that the Grunsky operator is Hilbert-Schmidt. In particular, we show that if the Grunsky operator is Hilbert-Schmidt, then $\f(\disk)$ is a Weil-Petersson quasidisk. 
  The formulations of the results and proofs make essential use of a geometric treatment of Smirnov space as a space of half-order differentials.
\end{abstract}

\tableofcontents

\begin{section}{Introduction}
\begin{subsection}{Results and literature}

The Grunsky operator is an operator associated to a conformal map $\f$ from the unit disk into the complex plane (a biholomorphic function). It has been studied for many years in complex function theory, for example in association with extremal problems for univalent functions \cite{Duren_book,Pommerenkebook}, regularity of the boundary of the image domain \cite{Pommerenkebook,Shen}, and potential theory \cite{Schiffer_expository}. It can be formulated and viewed in many ways, for example as an integral operator \cite{BergmanSchiffer}, or in terms of generating functions \cite{Duren_book,Pommerenkebook}. It also appears in connection to Teichm\"uller theory and symplectic geometry \cite{TNT}.

The functional-analytic theory of the Grunsky operator has so far explicitly or implicitly involved the Dirichlet norm on functions, or equivalently, the Bergman norm on one-forms. In this paper, we define an analogue of the Grunsky operator on the Hardy space of the disk, and prove results for this new operator analogous to those for the standard Grunsky operator.
The interpretation of the Grunsky operator involves the Smirnov space on the image of $\f$ and its complement (which on the disk agrees with the Hardy space).  

The main results are as follows.
All of the results assume that $\f$ maps onto a domain with Ahlfors-regular rectifiable boundary $\Gamma$.
(1) The Grunsky operator is bounded; (2) the graph of the Grunsky operator is the pull-back of the Smirnov space of the complement (Thm.~\ref{thm:GraphOfGrunsky}); (3) if the Grunsky operator is Hilbert-Schmidt, then $\Gamma= \f(\mathbb{S}^1)$ is a Weil-Petersson class quasicircle (Cor.~\ref{cor:HSImpliesWP}).  (4) We also show a partial converse.

These results are largely  motivated by the Kirillov--Yuri'ev/Nag--Sullivan (or KYNS) period map of $\mathrm{Diff}(\mathbb{S}^1)/\text{M\"ob}(\mathbb{S}^1)$ into the infinite Siegel disk, and the representation theory of $\mathrm{Diff}(\mathbb{S}^1)$.  The infinite Siegel disk was conceived of and investigated by G. Segal \cite{SegalUnitary}, in association with representations on a symmetric Fock space. In that paper $\mathrm{Diff}(\mathbb{S}^1)$ acts symplectically on smooth functions on the circle by composition, thus embedding $\mathrm{Diff}(\mathbb{S}^1)$ into the Lagrangian Grassmannian.  The connection with the Grunsky operator comes from conformal welding: an element $\phi$ of $\mathrm{Diff}(\mathbb{S}^1)$ (or more generally, a quasisymmetry of $\mathbb{S}^1$) can be written $\phi=\g^{-1} \circ \f$ where $\f$ and $\g$ are conformal maps of the disk and its complement in the sphere respectively.  A. Kirillov and D. Yuri'ev \cite{KY2} showed that the Lagrangian corresponding to $\phi$ is the graph of the Grunsky matrix of $\f$. It was shown by S. Nag and D. Sullivan \cite{NagSullivan} that the symplectic action by composition extends to quasisymmetries of $\mathbb{S}^1$, where smooth functions are replaced by the homogeneous Sobolev space $\dot{H}^{1/2}(\mathbb{S}^1)$; indeed quasisymmetries are precisely the bounded symplectomorphisms. In particular, the KYNS period map embeds the universal Teichm\"uller space into the symplectic Lagrangian Grassmannian.
The second two authors proved in \cite{Schippers_Staubach_Grunsky_expository} that when $\f(\disk)$ is a quasidisk, the graph of the Grunsky operator is the pull-back under $\f$ of the homogeneous Dirichlet space of the complement of $\f(\disk)$. It was furthermore shown by Y.~Shen \cite{Shen} and L.~Takhtajan and L--P.~Teo \cite{TNT} independently that the Grunsky operator is Hilbert--Schmidt if and only if the conformal map $\f$ is in what is called the Weil--Petersson class.
The Hilbert--Schmidt condition means that the corresponding Lagrangian is in the restricted Lagrangian Grassmannian, which has importance in the representation theory. See \cite[Section 6]{WP_Thurston_Schippers_Staubach} for a survey in the context of Weil--Petersson Teichm\"uller theory.   

The action of $\mathrm{Diff}(\mathbb{S}^1)$ on function spaces on the circle is an important ingredient in Segal's formulation of conformal field theory, and representation theory of the Virasoro algebra \cite{Huang,KW22,Segal04,SegalUnitary,Tener17}. The actions in the previous paragraph are symplectic, and not orthogonal.  So the cited results do not apply to the orthogonal Lagrangian Grassmannians of $L^2(\mathbb{S}^1)$ which play a role in fermionic models, see e.g. Segal \cite{Segal04}, J. Tener \cite{Tener17}, or P. Kristel and C. Waldorf \cite{KW22}.
It is then natural to ask whether there are fermionic analogues of the original Grunsky operator --- which we could call bosonic --- and corresponding results about its geometric and algebraic meaning.
It turns out that although an appropriate analogue of the closely related Faber operator (which is a composition of certain pull-back and a certain Cauchy-type integral operator)  existed in the literature, an analogue of the Grunsky operator did not. We advance such an analogue, and show that its graph is the Grassmannian associated to the conformal map $\f$, see \cref{thm:GraphOfGrunsky}.  We also show that if the Grunsky operator is Hilbert-Schmidt (among conformal maps onto domains with rectifiable and Ahlfors-regular boundary), then the conformal map is in the Weil-Petersson class (\cref{cor:HSImpliesWP}). As in the symplectic case, this condition has importance in representation theory and fermionic conformal field theory. 

We adopt a presentation which can be understood entirely in terms of complex analytic function theory (see below). However, for those who wish to compare it with the literature in conformal field theory, which in the fermionic theories involves spin bundles over domains and curves in the Riemann sphere \cite{Tener17,KW22}, we have provided \cref{sec:HalfOrderGeo}. It can also be asked whether there are generalizations to Riemann surfaces of the fermionic Grunsky operator, analogous to the generalization of the bosonic Grunsky operator obtained by M. Shirazi \cite{Shirazi_Grunsky}. Such a generalization would have to be formulated in terms of spin bundles over Riemann surfaces, either implicitly or explicitly.

Aside from these results, a major part of this paper is devoted to demonstrating that the operator we define here is indeed the natural analogue of the original Grunsky operator that was associated to one-forms (the so-called bosonic Grunsky operator). This can be done even from a complex function-theoretic point of view, partly thanks to existing analogues of the Faber operator.  While the classical Grunsky operator can be viewed as acting on Bergman/Dirichlet spaces of one-forms/functions, this new Grunsky operator acts on the Smirnov space of half-order differentials.

To accomplish this, we make use of a perspective of M. Bolt and D. Barrett on the Cauchy and related operators in terms of these half-order differentials, which they used to investigate the geometry of the Kerzman-Stein operator \cite{Barrett_Bolt_Mobius,Barrett_Bolt_Laguerre}. The point of view seems to us natural both geometrically and analytically. The half-order differential point of view was also the centre of a lengthy investigation by N. Hawley and M. Schiffer \cite{Hawley_Schiffer_half-order} into domain functions associated to conformal maps, invariant functions on Riemann surfaces, and the Schwarzian derivative and the Ricatti equation, among other things. In order to demonstrate that our Grunsky operator on Hardy space is the natural analogue of the classical Grunsky operator on Dirichlet/Bergman space, we also redevelop known results for the Faber operator and series from the half-order differential point of view of Barrett and Bolt.  Here we make no claims to originality except in the manner of presentation of known results, and also in the introduction of the ``overfare'' operator into the formulation of the Faber operator.  An overfare operator was used by the second two authors in the Dirichlet/Bergman space in investigation of the jump decomposition, Faber operator, and Grunsky operator for domains bounded by quasicircles \cite{Schippers_Staubach_Grunsky_expository}. Taken as a whole, this motivates the definition of our new Grunsky operator in terms of the Faber operator, in a way which is entirely analogous to the classical Grunsky operator. 

Although we make no use of the interesting theorems in \cite{Barrett_Bolt_Mobius,Barrett_Bolt_Laguerre,Hawley_Schiffer_half-order}, the approach plays a major role in this paper. We include a development of the formalism in the present context, since it motivates and clarifies the definitions and results -- indeed it led us to their formulation. Furthermore, in our opinion, it has considerable explanatory power and value on its own and provokes many new questions.  \\

{{\bf{Acknowledgements.}} The first and second authors were partially supported by the National Sciences and Engineering Research Council of Canada.
The first author gratefully acknowledges support from the Hausdorff Center for Mathematics.
 The third author is grateful to Andreas Str\"ombergsson for partial financial support through  a grant from Knut and Alice Wallenberg Foundation.
}

\end{subsection}
\begin{subsection}{Preliminaries}

Here we  recall some basic facts from the theory of quasiconformal Teichm\"uller spaces, theory of function spaces on the disk and planar domains, the theory of Cauchy integrals on Ahlfors regular curves (or rather its ramifications),  and finally weighted norm inequalities for singular integral operators, that will all be used in various sections of this paper.\\ 

{\bf{Notations and nomenclatures.}}\\
In this paper, a ``conformal'' map is a biholomorphism from a domain onto its image. Also, if the values of constants $C$ in estimates of the form $a\leq Cb$ are of no significance for our main purpose, then we use the notation $a\lesssim b$ as a shorthand for $a\leq Cb$. 
\begin{subsubsection}{Teichm\"uller spaces and the Weil-Petersson class}\label{Teichsection}

    Let $\overline{\mathbb{C}}$ denote the Riemann sphere, $\disk := \{ z \in \C \, \colon \, |z| < 1\}$ the unit disk, and \(\disk^{*}=\overline{\mathbb{C}}\setminus\mathrm{cl}({\disk})\) be the exterior of the unit disk, where $"\mathrm{cl}"$ denotes the closure.

Now, let \(L^{\infty}(\disk^{*})_1\) denote the open unit ball of the Banach space \(L^{\infty}(\disk^{*})\) of bounded
measurable functions on \(\disk^*\). For \(\mu \in L^{\infty}(\disk^*)_1\), 
extend it to \(\mathbb{D}\) by the reflection
\begin{equation}
\mu(z)=\overline{\mu\left(\frac{1}{\bar{z}}\right)} \frac{z^{2}}{\bar{z}^{2}},\qquad z \in \mathbb{D}.
\end{equation}
Let \(w_{\mu}\) be the unique quasiconformal homeomorphism from $\mathbb{C} \to \mathbb{C}$ (i.e.~\(\overline{\partial} w_{\mu}= \mu\, \partial w_{\mu}\) in $\mathbb{C}$) which fixes the points $-1$, $i$ and $1$. One defines an equivalence relation on the space of complex dilatations in $L^{\infty}(\disk^*)_1$ in the following way: \(\mu\) and \(\nu\) are equivalent if \(\left. w_{\mu}\right|_{\mathbb{S}^1} = \left. w_{\nu}\right|_{\mathbb{S}^1} \).

\begin{definition}\label{Teich}
The \emph{universal Teichm\"uller space}  \(T(1)\) is defined as 
\begin{equation}
    T(1)= \{[\mu] \colon \, \mu\in L^{\infty}(\disk^*)_1  \},
\end{equation}
 where \([\mu]\) denotes
the equivalence classes of $\mu$ according to the equaivalence relation above.
\end{definition}
Another model for the universal Teichm\"uller space is given as follows. First, we extend $\mu\in L^{\infty}(\disk^*)_1$  to be zero outside $\disk^*$. Then one considers the unique quasiconformal mapping $w^\mu$ that is the solution of the Beltrami equation $\overline{\partial} w^{\mu}= \mu\, \partial w^{\mu}$ and normalized by the conditions
$(w^{\mu}|_{\disk})(0)=0$, $(w^{\mu}|_{\disk})'(0)=1$, $(w^{\mu}|_{\disk})''(0)=0$.
Now define an equivalence relation $\sim$ on the space of complex dilatations in $L^{\infty}(\disk^*)_1$ in the following way: $\mu\sim\nu$ if $w^{\mu}\left|_{\disk}=w^{\nu}\right|_ {\disk}.$ Using this equivalence relation, one can define $T(1)=L^{\infty}(\disk^*)_1/\sim$. However we observe that since $w^{\mu}|_\disk= w^{\nu}|_\disk \iff w_{\mu}|_{\mathbb{S}^1}= w_{\nu}|_{\mathbb{S}^1}$, the two definitions of $T(1)$ are equivalent to each other.\\

These two models can be used to define a conformal welding associated to an element $[\mu]\in T(1)$ (which will be used in Section \ref{WP}).
Indeed, if $[\mu] \in T(1)$, then there is a corresponding conformal welding, which for appropriate quasiconformal mappings $\mathrm{g}_{ \mu}$ and $\mathrm{f}^{\mu}$ (on $\mathbb{C}$), is given by 
\begin{equation}\label{confweld}
    w_{\mu}=\mathrm{g}_{ \mu}^{-1} \circ \mathrm{f}^{\mu},
\end{equation}
 see \cite[Section 2.2]{TNT} and \cite{Lehto} for details.
We also note that $\mathrm{f}^{\mu}$ is conformal inside the unit disk and $\mathrm{g}_\mu$ is conformal outside the unit disk.

The universal Teichm\"uller space \(T(1)\) is a group (not a topological group)
under the composition of the quasiconformal mappings. The group law on
\(L^{\infty}\left(\mathbb{D}^{*}\right)_{1}\) is defined implicitly by

$$\lambda=\nu\ast \mu^{-1},$$
through $w_\lambda= w_\nu \circ w_{\mu}^{-1}$, where $\mu^{-1}$ is defined by the property that $\mu\ast \mu^{-1}=0$. Explicitly, the group law is given by
\begin{equation}
\lambda=\left(\frac{\nu-\mu}{1-\bar{\mu} \nu} \frac{\left(\partial w_{\mu}\right)}{\left(\overline\partial\bar{w}_{\mu}\right)}\right) \circ w_{\mu}^{-1}.
\end{equation}

If $\Phi: L^{\infty}(\disk^*)_1\to T(1)$ is the natural projection-map $\mu\mapsto [\mu]$, then the group structure on $L^{\infty}(\disk^*)_1$
projects to $T(1)$ by $[\lambda]\ast [\mu]=[\lambda \ast \mu]$ and the right translation defined by
\begin{equation}\label{defn:righttrans}
R_{[\mu]}: T(1) \longrightarrow T(1), \quad[\lambda] \longmapsto[\lambda * \mu],
\end{equation}
is a biholomorphic automorphism of $T(1).$ Moreover, for \(\mu \in L^{\infty}\left(\mathbb{D}^{*}\right)_{1}\) the tangent space at $[\mu]$ of $T(1)$ is given by
\begin{equation}
T_{[\mu]} T(1)=D_{0} R_{[\mu]}\left(T_{0} T(1)\right) 
\end{equation}

As was shown by L. Takhtajan and L-P. Teo \cite{TNT} $T(1)$ possesses a Hilbert manifold structure with a natural Hermitian metric. Namely, if
\begin{equation}\rho(z)=\frac{4}{\left(1-|z|^{2}\right)^{2}}
\end{equation}
denotes the the density of the hyperbolic metric on $\mathbb{D}^*$, and $dA_z$ denotes the Lebesgue area-measure on $\mathbb{C},$ define the Hilbert space of harmonic Beltrami differentials on $\mathbb{D}^*$ by 
\begin{equation}
H^{-1,1}\left(\mathbb{D}^{*}\right)=\left\{\mu=\frac{ \overline{\phi(z)}}{\rho(z)}, \phi \text { holomorphic on } \mathbb{D}^{*};\, \Vert \mu\Vert_{2} <\infty\right\},
\end{equation}
where
\begin{equation}\label{norm of mu}
   \Vert \mu\Vert^2_{2}:= \iint_{\mathbb{D}^{*}}|\mu|^{2} \rho(z)\, dA_z. 
\end{equation}

Then it was shown by Takhtajan and Teo \cite{TNT}
that $T(1)$ is a Hilbert manifold with uncountably many connected components. The tangent spaces are given by 
\begin{equation}
T(1) \ni[\mu] \mapsto D_{0} R_{[\mu]}\left(H^{-1,1}\left(\mathbb{D}^{*}\right)\right) \subset T_{[\mu]} T(1),
\end{equation}
where above $T_{[\mu]}T(1)$ denotes the tangent space at $[\mu]$ with respect to the classical Banach manifold structure.

Let $T_0(1)$ be the component of origin of the Hilbert manifold $T(1)$. As was shown in \cite{TNT}, $T_0 (1)$ is a subgroup of $T(1)$, it is a Hilbert manifold and, as opposed to $T(1)$, a topological group. Given $[\mu]\in T(1)$ with the corresponding conformal welding \eqref{confweld}, it was shown in \cite{TNT} that $[\mu]\in T_0(1)$ if and only if 
\begin{equation}\label{Schwarzian WP condition}
   \int_{\disk} \frac{|\mathcal{S}\mathrm{f}^{\mu}(z)|^2}{\rho(z)}\, dA_z <\infty, 
\end{equation}
where $\mathcal{S}f:= \left(f^{\prime \prime} / f^{\prime}\right)^{\prime}-1 / 2\left(f^{\prime \prime} / f^{\prime}\right)^{2}$ is the Schwarzian derivative of a function $f$. This leads to the following definition.
\begin{definition}\label{defn:WP}
    The class of such $[\mu]$'s in $T_0(1)$ is called the \emph{Weil-Petersson class}. Note that in this case $\Vert \mu\Vert_2 <\infty.$ With a slight abuse of nomenclature, we also say that $\f^\mu$ belongs to the Weil-Petersson class (referred to as the WP--class). In this case the quasicircle $\f^\mu(\mathbb{S}^1)$ is referred to as a WP--class quasicircle.
\end{definition}

\end{subsubsection}

\begin{subsubsection}{Function spaces on the disk}
     Next, we recall some facts about certain function spaces on the disk and planar domains that are used in this paper.

 \begin{definition}\label{defn:Bergmanspace}
The \emph{Bergman space} of the disk, denoted $\mathcal{A}(\disk)$, is the Hilbert space consisting of holomorphic functions $f: \disk \rightarrow \C$ such that 
\begin{equation}
    \Vert f\Vert^2_{\mathcal{A}(\disk)}:= \int_{\disk} |f(z)|^{2}\, dA_z <\infty.
   \end{equation} 
\end{definition}    
    
\begin{definition}\label{defn:Hardyspace}
For $0<p<\infty$, the \emph{Hardy spaces} of the disk, denoted $H^{p}(\disk)$, are the Banach spaces consisting of holomorphic functions $h: \disk \rightarrow \C$ such that 
   \begin{equation}\label{eq:HardyNorm}
     \Vert h\Vert^p_{H^p(\disk)}:=  \sup_{0<r<1} \int_0^{2 \pi} |h(r e^{i\theta})|^{p} d\theta <\infty.
   \end{equation}
   
\end{definition}
It is a well-known fact of function theory that
\begin{equation}\label{cont inclusion hardy-bergman}
    H^2(\disk)\hookrightarrow \mathcal{A}(\disk),
\end{equation}
where $\hookrightarrow$ denotes continuous (i.e.~bounded) inclusion.\\
Let $\Omega$ be a simply-connected domain in the plane, which is conformally equivalent to the unit disk.
 Let $g_q$ denote its Green's function with singularity at $q \in \Omega$.
 For $0<r<1$, let $\Gamma_{q,r}$ denote the level curve \begin{equation}\label{greenlevel}
    \Gamma_{q,r} = \{ z \in \Omega : g_q(z) = - \log{r} \}. 
 \end{equation}   
 \begin{definition}
    Given a rectifiable Jordan curve $\Gamma$ with complementary component $\Omega$,  one defines the \emph{generalized Hardy space} $\hat{H}^2(\Omega)$ to be closed subspace of $L^2(\Gamma,|dz|)$ such that 
  \[ \int_{\Gamma} f(z) z^k dz =0 , \ \ \ k \geq 0 \]  
  in the case that $\Omega$ is the bounded component of the complement, and
  \[ \int_{\Gamma} f(z) z^k dz =0 , \ \ \ k \leq -1 \]  
  in the case that $\Omega$ is the unbounded component.  
 \end{definition}

   A class of functions that play an important role in approximation theory is the so-called Smirnov class.
\begin{definition}
   Let $0<p<\infty$. Let $\Omega$ be a simply-connected domain conformally equivalent to the disk, and let $h: \Omega \rightarrow \C$ be holomorphic.
   We say that $h \in E^p(\Omega)$ (the \emph{Smirnov space}), if there is a sequence of rectifiable simple closed curves $\Gamma_n$ eventually enclosing any compact subset of $\Omega$ such that 
   \begin{equation}\label{eq:SmirnovSup}
    \Vert f\Vert^p_{E^p(\disk)}:=    \sup_n \int_{\Gamma_n} |f(z)|^p |dz| <\infty.
    \end{equation}
\end{definition}

  If this holds for some $\Gamma_n$, it also holds in particular for $\Gamma_{q,r}$ for any $q$, see \cite[Theorem 10.1]{Duren_Hp_book}.\\

A \emph{Smirnov domain} is a bounded simply-connected domain \(\Omega\) with a rectifiable Jordan boundary in the complex plane \(\mathbb{C}\) with the following property: there is a conformal mapping \(z=\phi(w)\) from the disk \(|w|<1\) onto \(\Omega\) such that for \(|w|<1\) the harmonic function \(\log \left|\phi^{\prime}(w)\right|\) can be written as the Poisson integral of its non-tangential boundary values \(\log \left|\phi^{\prime}\left(e^{i \theta}\right)\right|:\)

\begin{equation}
\log \left|\phi^{\prime}\left(r e^{i \theta}\right)\right|=\frac{1}{2 \pi} \int^{2 \pi} \frac{1-r^{2}}{1+r^{2}-2 r \cos (t-\theta)} \log \left|\phi^{\prime}\left(e^{i t}\right)\right| d t
\end{equation}
  
  If $\Omega$ is a Smirnov domain, then $\hat{H}^2(\Omega) = E^2(\Omega)$, see \cite[Chapter 12]{CoifmanMeyer}.  
 \end{subsubsection}
 
\begin{subsubsection}{Ahlfors regular curves and G. David's theorem} 
In this paper, the regularity conditions that are required of the boundary curves are rectifiability and Ahlfors regularity. The latter is defined as follows.
\begin{definition}
    Let \(\Gamma\) be a Borel set in \(\mathbb{R}^{2}\). We say that \(\Gamma\) is a
\emph{Ahlfors-regular} if it is bounded and if there is a constant \(C_{\Gamma}\) such
that
\begin{equation}
\frac{r}{C_{\Gamma}}  \leq \mathscr{H}^{1}(B(x, r) \cap \Gamma) \leq C_{\Gamma} r
\end{equation}
for all \(x \in \Gamma, 0<r \leq 1\), where \(\mathscr{H}^{1}\) denotes the \(1\)-dimensional Hausdorff measure.
\end{definition} 
For a Jordan curve $\Gamma$ which splits the plane into complementary components, $\Omega_1$ and $\Omega_2$, A.~Calder\'on \cite{Calderon} posed the problem of whether the space $L^2(\Gamma)$ is the direct sum of the generalized Hardy spaces  $\hat{H}^2(\Omega_1)$ and $\hat{H}^2(\Omega_2)$, when these spaces are realized as subspaces of $L^2(\Gamma)$ via their corresponding trace operators.
G.~David \cite{David} characterized the curves $\Gamma$ for which Calder\'on's problem has an affirmative answer. David's result is
 \begin{theorem}\label{Davids} Let $\Gamma$ be a rectifiable curve in the plane, of finite total length, and let $\Omega_1$ and $\Omega_2$ be its complementary components in the Riemann sphere. 
 The direct sum decomposition 
 \[  L^2(\Gamma,|dz|) = \hat{H}^2(\Omega_1) \oplus \hat{H}^2(\Omega_2)   \]
 holds if and only if $\Gamma$ is Ahlfors-regular.  The decomposition is the jump decomposition obtained from the Cauchy integral. 
 \end{theorem}
 
\end{subsubsection}
\begin{subsubsection}{Muckenhoupt weights and weighted norm inequalities}\label{Muckweights}
    We recall very briefly some basic facts about weighted norm inequalities, that will be used in Section \ref{WP}.
    \begin{definition}
        For a fixed \(1<p<\infty\), one says that a non-negative function \(w: \mathbb{R}^{n} \rightarrow[0, \infty)\) belongs to the \emph{Muckenhoupt} \(A_{p}\)-\emph{class}, if \(w\) is locally integrable and there is a constant \(C\)
such that, for all balls \(B\) in \(\mathbb{R}^{n}\), one has
$$
\left(\frac{1}{|B|} \int_{B} w(x) d x\right)\left(\frac{1}{|B|} \int_{B} w(x)^{-\frac{p'}{p}} d x\right)^{\frac{p}{p'}} \leq C<\infty
$$
where \(|B|\) is the Lebesgue measure of \(B\), and \(p'\) is the H\"older conjugate of $p$, meaning that \(\frac{1}{p}+\frac{1}{p'}=1.\)
    \end{definition}
    It is well-known that for $p>1$ and $x\in\mathbb{R}^n$, the weights $|x|^{\alpha} \in A_{p}$ if and only if $-n<\alpha<n(p-1).$\\
    
    One says that the function $f\in L^{p}_{w}$, if $f$ is measurable and $$\Vert f\Vert_{L^p_w}:= \Big\{\int_{\mathbb{R}^n} |f(x)|^p \, w(x) \, dx\Big\}^{1/p}<\infty.$$
    Apart from the significant role that the Muckenhoupt weights play in the weighted $L^p$--boundedness of Hardy-Littlewood maximal operators, it was shown by R. Coifman and C. Fefferman \cite{Coifman-Fefferman} that singular integral operators of Calder\'on-Zygmund type are also bounded on weighted $L^p$ spaces equipped with Muckenhoupt weights. More precisely
    \begin{theorem}\label{CF theorem}
        Let $T$ be a Calder\'on-Zygmund operator of convolution-type. Then for $1<p<\infty$ and $w\in A_p$, one has the weighted norm inequality
$$\Vert Tf\Vert_{L^p_w}\leq C_p \Vert f\Vert_{L^p_w}.$$
    \end{theorem}
    
    Note that the Hilbert and the Beurling transforms
    \begin{align}\label{Beurling}
        \mathscr{H} f(x)&= \textnormal{P.V.} \int_{\mathbb{R}} \frac{f(y)}{x-y} dy, & &\text{ and } & \mathscr{B} f(z)&= \textnormal{P.V.}\int_{\mathbb{C}}\frac{f(\zeta)}{(\zeta-z)^2} \, dA_u
    \end{align}
    respectively, are examples of Calder\' on-Zygmund operators to which this result applies.
\end{subsubsection}

\end{subsection}
\end{section}
\begin{section}{Smirnov spaces of half-order differentials}
\begin{subsection}{Half-order differentials}\label{sec:HalfOrderDifferentials}

    We give a nuts-and-bolts description of holomorphic half-order differentials, that is in line with a common function-theoretic way of thinking of differentials (e.g.~\cite{Lehto}).
    In Appendix \ref{sec:HalfOrderGeo}, we show that half-order differentials are really sections of the square-root of the canonical line bundle, and moreover, that all definitions in this section are completely compatible with this point of view.
    While this differential-geometric viewpoint is not strictly necessary to understand or prove our results, it is a useful source of intuition, and moreover provides our results with a broader context.

    We shall define holomorphic $\frac{1}{2}$-differentials in terms of how they transform under holomorphic change of coordinates.
    Whereas the transformation law of a 1-differential involves the derivative of the holomorphic transformation, for $\frac{1}{2}$-differentials we require the square root of the derivative.
    First, we make precise what we mean by this.
    Recall that if $f: \Omega_{2} \rightarrow \Omega_{1}$ is a biholomorphism of simply-connected (proper) domains, then there exist exactly two holomorphic functions $g: \Omega_{2} \rightarrow \C^{\times}$ that satisfy $g(z)^{2} = f'(z)$ for all $z \in \Omega_{2}$ (here $\C^{\times}= \C\setminus \{0\}$).
    We denote by $\mathrm{Conf}(\Omega_{2},\Omega_{1})$ the set of biholomorphisms from $\Omega_{2}$ to $\Omega_{1}$ and define
    \begin{equation*}
        \Hdiffhat(\Omega_{2},\Omega_{1}) := \{ (f,g) \colon \, f \in \mathrm{Conf}(\Omega_{2},\Omega_{1}), \, g(z)^{2} = f'(z) \}.
    \end{equation*}
    We will denote elements $\hat{f} \in \Hdiffhat(\Omega_2,\Omega_1)$ by $\hat{f} = (f,\sqrt{f'})$ where $f$ is a biholomorphism from $\Omega_{2}$ to $\Omega_{1}$, and it is understood that $\sqrt{f'}$ denotes a definite choice of one of the two possible branches of square root of $f'$.

    \begin{definition}\label{def:HalfOrderDifferentials}
        The space of \emph{holomorphic half-order differentials} on $\Omega_{1}$, denoted $\Omega^{\frac{1}{2},0}(\Omega_{1})$ is the space of holomorphic functions from $\Omega_{1}$ to $\C$.
        Elements of $\Omega^{\frac{1}{2},0}(\Omega_{1})$ are denoted by $h \sqrt{dz} = h dz^{1/2}$, where $h: \Omega_{1} \rightarrow \C$ is a holomorphic function.
        If $\hat{f} = (f,\sqrt{f'}) \in \Hdiffhat(\Omega_{2},\Omega_{1})$, then for $h \sqrt{dz} \in \Omega^{\frac{1}{2},0}(\Omega_{1})$, we define
        \begin{equation}\label{eq:HalfOrderTransformation}
            \hat{f}^* (h \sqrt{dz}) \defeq (h \circ f) \sqrt{f'} \sqrt{dz} \in \Omega^{\frac{1}{2},0}(\Omega_{2}),
        \end{equation}
        or, equivalently
        \begin{equation*}
            \hat{f}^* (h \sqrt{dz}) (z) = h(f(z)) \sqrt{f'(z)} \sqrt{dz}, \qquad z \in \Omega_{2}.
        \end{equation*}
        We equip $\Omega^{\frac{1}{2},0}(\Omega)$ with the bilinear pairing
        \begin{equation}\label{eq:HalfOrderPairing}
            \Omega^{\frac{1}{2},0}(\Omega) \times \Omega^{\frac{1}{2},0}(\Omega) \rightarrow \Omega^{1,0}(\Omega), \quad
            (h_{1} \sqrt{dz}, h_{2} \sqrt{dz}) \mapsto h_{1}h_{2} dz.
        \end{equation}
    \end{definition}
    The symbol $\sqrt{dz}$ is formal, but it serves to remind us of Equations \eqref{eq:HalfOrderTransformation} and \eqref{eq:HalfOrderPairing}.
    The space $\Hdiffhat(\Omega) \defeq \Hdiffhat(\Omega,\Omega)$ is a group, and in fact it is a double cover of the group of biholomorphisms from $\Omega$ to itself.
    Equation \eqref{eq:HalfOrderTransformation} then yields a group action of $\Hdiffhat(\Omega)$ on $\Omega^{\frac{1}{2},0}(\Omega)$.

    \begin{definition}
        The space of \emph{anti-holomorphic half-order differentials} on $\Omega_{1}$, denoted $\Omega^{0,\frac{1}{2}}(\Omega_{1})$ is the space of anti-holomorphic functions from $\Omega_{1}$ to $\C$.
        Elements of $\Omega^{0,\frac{1}{2}}(\Omega_{1})$ are denoted by $h \sqrt{d\overline{z}}$.
        If $\hat{f} = (f,\sqrt{f'}) \in \Hdiffhat(\Omega_{2},\Omega_{1})$, then for $h\sqrt{d\overline{z}} \in \Omega^{0,\frac{1}{2}}(\Omega_{1})$, one defines
        \begin{equation*}
            \hat{f}^*(h \sqrt{d\overline{z}}) \defeq (h \circ f) \overline{\sqrt{f'}} \sqrt{d\overline{z}} \in \Omega^{0,\frac{1}{2}}(\Omega_{2}).
        \end{equation*}
    \end{definition}

    At the moment, the space $\Omega^{\frac{1}{2},0}(\Omega)$ is simply a vector space.
    In the sequel, we consider a subspace consisting of elements satisfying a certain integrability condition, see Eq.~\eqref{eq:norm_definition}, and see that this leads naturally to Smirnov spaces.

\end{subsection}

\begin{subsection}{The Smirnov space of half-order differentials}

In this section we add some regularity to the half-order differentials (introduced in Sec.~\ref{sec:HalfOrderDifferentials}) to obtain a model of the Smirnov space. In brief, we view the Smirnov space as half-order differentials of the form $h(z) \sqrt{dz}$ where $h(z) \in E^2(\Omega)$.  We recall some basic results, in order to establish our terminology.

  Let $\mathcal{A}^{1/2}(\Omega,q)\subset \Omega^{\frac{1}{2},0}(\Omega)$ denote the space of holomorphic half-order differentials $h\, \sqrt{dz}$ satisfying 
  \begin{equation} \label{eq:norm_definition}
    \| h \sqrt{dz} \|_{\Omega,q}^2 := \lim_{r \nearrow 1} \int_{\Gamma_{q,r}} |h(z)|^2|dz| <\infty,    \end{equation}
    where $\Gamma_{q,r}$ is defined in \eqref{greenlevel}.
  Equivalently, using the notation in \eqref{greenlevel}, the conformal invariance of Green's function and the fact that Green's function of the disk satisfies $g_0(z)=-\log|z|$, one may write the curves $\Gamma_{q,r}$ as curves $f(|z|=r)$ where $f:\disk \rightarrow \Omega$ is a conformal map such that $f(0)=q$.  This yields that 
  \begin{equation} \label{eq:norm_definition_conformal_map}
    \| h \, \sqrt{dz} \|_{\Omega,q}^2 := \lim_{r \nearrow 1} \int_{f(|z|=r)} |h(z)|^2 \,|dz|.    
  \end{equation}
 
  Thus, with \eqref{eq:norm_definition} and this definition, we can identify $\mathcal{A}^{1/2}(\Omega,q)$ with the Smirnov space $E^2(\Omega)$
  \[   \mathcal{A}^{1/2}(\Omega,q) = \left\{ h \, dz^{1/2} : h \in E^2(\Omega) \right\}.   \]
  Furthermore, as a collection of functions, $\mc{A}^{1/2}(\Omega,q)$ is independent of $q$.
  We shall shortly show that the norm does not depend on $q$ either, that is, we shall prove the following result.
    \begin{proposition}\label{prop:HalfOrderNorm}
        Let $\Omega$ be a simply-connected domain.
        Then, we have, for all $h \, \sqrt{dz} \in \Omega^{\frac{1}{2},0}(\Omega)$
        \begin{equation*}
            \| h \sqrt{dz} \|_{\Omega,p} = \| h \sqrt{dz} \|_{\Omega,q}
        \end{equation*}
        for all $p,q \in \Omega$.
    \end{proposition}
    
    We will prove Prop.~\ref{prop:HalfOrderNorm} below, but first we take some preliminary steps.
    \begin{lemma}\label{lem:Isometries}
        If $\hat{f} \in \Hdiffhat(\Omega_{2},\Omega_{1})$, then the map $\hat{f}^{*}: \Omega^{\frac{1}{2},0}(\Omega_{1}) \rightarrow \Omega^{\frac{1}{2},0}(\Omega_{2})$ restricts to an isometry
        \begin{equation*}
        {\hat{f}^*}: \mc{A}^{1/2}(\Omega_{1},p) \rightarrow \mc{A}^{1/2}(\Omega_{2},f^{-1}(p)),
        \end{equation*}
        for any $p \in \Omega_{1}$.
    \end{lemma}
    \begin{proof}
        Let $h \sqrt{dz} \in \mc{A}^{1/2}(\Omega_{1},p)$, and $\hat{f} = (f, \sqrt{f'}) \in \Hdiffhat(\Omega_{2},\Omega_{1})$.
        The crucial observation here is that we have, for all $r \in (0,1)$ that $f(\Gamma_{f^{-1}(p),r}) = \Gamma_{p,r}$. This yields for $r \in (0,1)$ that
        \begin{equation*}
            \int_{\Gamma_{f^{-1}(p),r}}| h(f(z)) | |f'(z)| |dz| = \int_{\Gamma_{p,r}}| h(w) | |dw|.
        \end{equation*}
        Taking the limit as $r \nearrow 1$ on both sides, yields the desired conclusion.
    \end{proof}
    
   We note that, by comparing Eq.~\eqref{eq:SmirnovSup} with Eq.~\eqref{eq:HardyNorm}, we obtain an identification of Banach spaces
   \begin{equation}\label{Hardy_Ahalf correspondence}
       H^{2}(\disk) = \mc{A}^{1/2}(\disk,0).
   \end{equation}
   This identification will allow us to exploit some well-known results in the context of Hardy spaces to prove \cref{prop:HalfOrderNorm}.
   
    \begin{proof}[Proof of \textnormal{Proposition \ref{prop:HalfOrderNorm}}]
    We first consider the case that $\Omega = \disk$.
    Let $h \sqrt{dz} \in \mc{A}^{1/2}(\disk,0)$ be arbitrary, and let
  $C_r$ be the circle $|z|=r$ for $0<r<1$.
  As is well-known, elements of $H^2(\disk)$ have non-tangential limits almost everywhere and by the mean convergence theorem the Hardy space norm equals the $L^2$-norm of the boundary values with respect to the measure $d\theta$.
  By the (isometric) identification of $\mc{A}^{1/2}(\disk,0)$ with $H^{2}(\disk)$, we thus have
  \begin{equation*}
      \| h\, \sqrt{dz} \|^2_{\disk,0} = \lim_{r \nearrow 1} \int_{C_r} |h(z)|^2 |dz| = \int_{\mathbb{S}^1} |h(z)|^2 |dz|.
  \end{equation*}
  For an arbitary M\"obius transformation $M$ preserving $\disk$, we then have 
  \begin{align} \label{eq:Hardy_p_independent}
     \int_{\mathbb{S}^1} |h(z)|^2 |dz| &= \int_{\mathbb{S}^1} |h(M(z))|^2|M'(z)| |dz| \nonumber \\ & = \lim_{r \nearrow 1} \int_{C_r} |h(M(z))|^2|M'(z)| |dz| = \lim_{r \nearrow 1} \int_{M^{-1}(C_r)} |h(z)|^2| |dz| \nonumber \\
     & = \| h\, \sqrt{dz} \|^2_{\disk,M^{-1}(p)}.  
   \end{align}
    Now, let $\Omega$ be an arbitrary simply-connected domain, and let $p,q \in \Omega$.
    Let $\hat{f} = (f,\sqrt{f'}) \in \Hdiffhat(\disk, \Omega)$ be such that $f^{-1}(p) = 0$.
    We then have, using \cref{lem:Isometries},
    \begin{equation*}
        \| h \sqrt{dz} \|^{2}_{\Omega,p} = \| \hat{f}^* h \sqrt{dz} \|^{2}_{\disk, 0} = \| \hat{f}^*  h \sqrt{dz} \|^{2}_{\disk, f^{-1}(q)} = \| h \sqrt{dz} \|^{2}_{\Omega,q}. \qedhere
    \end{equation*}
   \end{proof}

  \begin{remark}
   In many sources it is required that the boundary of $\Omega$ be rectifiable in order to define the Smirnov space. Here it is not necessary, but we will add that condition in the next section.
  \end{remark}

 Observe that the definition is entirely conformally invariant; the regularity of the boundary plays no role in either the norm or the space. 
  Observe also that the definition extends without problem to arbitrary simply-connected domains $\Omega$ in the Riemann sphere, so long as one observes that if $\infty \in \Omega$ then we must assume that for $h(z) \sqrt{dz} \in \mathcal{A}^{1/2}(\Omega)$ the function $h(1/z) /z$ is holomorphic at $0$.
  
  We also define an inner product on $\mathcal{A}^{1/2}(\Omega)$, namely 
  \[  \left( h_1\sqrt{dz},h_2 \sqrt{dz} \right)_\Omega = 
  \lim_{r \nearrow 1} \int_{\Gamma_{p,r}} h_1(z) \overline{h_2(z)} |dz|  \]
  which incidentally is nicely motivated by the product 
  \begin{equation} \label{eq:product}
   h_1(z) dz^{1/2} \cdot \overline{h_2(z)} d\bar{z}^{1/2} = h_1(z) \overline{h_2(z)} |dz|.  
  \end{equation}
  Arguing as above, we see that this is independent of $p$ (c.f.~\cref{prop:HalfOrderNorm}).
  Moreover, if $\hat{f} \in \Hdiffhat(\Omega_{2},\Omega_{1})$, then the associated map
  \begin{equation*}
      {\hat{f}^*}: \mc{A}^{1/2}(\Omega_{1}) \rightarrow \mc{A}^{1/2}(\Omega_{2})
  \end{equation*}
  is unitary (c.f.~\cref{lem:Isometries}).
  
  We denote by $\overline{\mathcal{A}^{1/2}(\Omega)}$ the set of differentials $\overline{h}\, d\bar{z}^{1/2}$ such that $h\, d{z}^{1/2} \in \mathcal{A}^{1/2}(\Omega)$, and define the inner product via  
  \[  \left( \overline{h_1} d\bar{z}^{1/2},\overline{h_2} d\bar{z}^{1/2} \right)_\Omega = 
  \lim_{r \nearrow 1} \int_{\Gamma_{p,r}} \overline{h_1(z)} {h_2(z)} |dz|.  \]
  \begin{definition}
  We set
   \[   \mathcal{A}^{1/2}_h(\Omega) :=  \mathcal{A}^{1/2}(\Omega) \oplus \overline{\mathcal{A}^{1/2}(\Omega)}, \]
  and extend the inner product so that the two subspaces are orthogonal.    
  \end{definition}

  \begin{remark} \label{re:heuristic_inner_product}
  If we also extend the product (\ref{eq:product}) in the obvious way, then we have that
  \begin{align*} \label{eq:general_pairing_on_sum}
   & \left( h_1 dz^{1/2} + \overline{H_1} d\bar{z}^{1/2}, h_2 dz^{1/2} + \overline{H_2} d\bar{z}^{1/2}  \right) \\
   &  \ \ \ \  = \left( h_1 dz^{1/2} ,h_2 dz^{1/2} \right) + 
    \left( \overline{H_1} d\bar{z}^{1/2}, \overline{H_2} d\bar{z}^{1/2} \right)\\ 
    & \ \ \ \  =  \lim_{r \nearrow 1} \int_{\Gamma_{p,r}} 
   \left( h_1(z) \overline{h_2(z)}  |dz| +  \overline{H_1(z)} H_2(z) |dz| + h_1(z) H_2(z) dz +   \overline{h_2(z)} \overline{H_1(z)} d\bar{z} \right) \\
   & \ \ \ \  =  \lim_{r \nearrow 0} \int_{\Gamma_{p,r}} \left(   h_1(z) dz^{1/2} + \overline{H_1(z)} d\bar{z}^{1/2}  \right) \left( \overline{h_2(z)} d\bar{z}^{1/2} + {H_2(z)} d {z}^{1/2}  \right), 
  \end{align*}
  because the third and the fourth integral in the third line vanish for all $0<r<1$ by holomorphicity of $h_k$ and $H_k$.
  Thus the choice that the holomorphic and anti-holomorphic  spaces are orthogonal is consistent with the product. In fact the product is the natural symmetric product arising in the differential geometric definition of differentials.   
  \end{remark}
  
    The obvious analog of \cref{lem:Isometries} for $\overline{\mc{A}^{1/2}(\Omega)}$ holds, from which it follows that elements of $\Hdiffhat(\Omega_{2},\Omega_{1})$ also yield isometries.
    In summary, we have the following proposition.
  \begin{proposition}
   \label{pr:pull-back_is_isometry}
   If $\hat{f} \in \Hdiffhat(\Omega_{2},\Omega_{1})$, then we have an isometry
   \begin{equation*}
      { \hat{f}^*}: \mc{A}_{h}^{1/2}(\Omega_{1}) \rightarrow \mc{A}_{h}^{1/2}(\Omega_{2}),
   \end{equation*}
   which sends $\mc{A}^{1/2}(\Omega_{1})$ to $\mc{A}^{1/2}(\Omega_{2})$ and $\overline{\mc{A}^{1/2}(\Omega_{1})}$ to $\overline{\mc{A}^{1/2}(\Omega_{2})}$.
  \end{proposition}

  \begin{remark}  
  Given an arbitrary simply-connected domain $\Omega$, one can define the Hardy space $H^2(\Omega)$ as the set of holomorphic functions on $\Omega$ with a harmonic majorant \cite{Duren_Hp_book}. For a conformal map $f:\Omega_1 \rightarrow \Omega_2$ 
  \[ h \in H^2(\Omega_2)  \Leftrightarrow  h \circ f \in H^2(\Omega_1), \]
  whereas a function $h$ is in the Smirnov space $E^2(\Omega_2)$ if and only if $h \circ f \sqrt{f'}$ is in $E^2(\Omega_1)$ \cite[10.1]{Duren_Hp_book}.   
  Thus elements of the Hardy space transform as functions while elements of the Smirnov space transform as half-order differentials. 
  \end{remark}
  \begin{remark}  \label{re:Smirnov_Hardy_confusion}
   Note that if the boundary of the domain is sufficiently regular, e.g.~$C^2$ smooth, then the Hardy spaces and Smirnov spaces agree  \cite[Theorem 10.2]{Duren_Hp_book}. Nevertheless, the distinction between their geometric natures should be kept in mind. 
  \end{remark}

\end{subsection}

 \begin{subsection}{Rectifiable curves and overfare} 
  Now assume that $\Omega$ is a simply connected domain whose boundary $\Gamma$ is rectifiable. In that case, the unit tangent $T$ exists almost everywhere on $\Gamma$, and so there is a well-defined notion of non-tangential limit almost everywhere on $\Gamma$. 
  
  We have the following well-known result (rephrased slightly).
  \begin{theorem} \label{th:nontangential_bvs_exist} Let $\Omega$ be a simply connected domain in the sphere whose boundary $\Gamma$ is rectifiable.
   Given any $h\, dz^{1/2} \in \mathcal{A}^{1/2}(\Omega)$, $h$ has a non-tangential limit almost everywhere on $\Gamma$.
   If the non-tangential limits vanish on a set of non-zero measure, then $h\, dz^{1/2}=0$.
   Finally,   
   \[    \int_{\Gamma} |h(z)|^{2}\, |dz| = \| h\, dz^{1/2} \|^2_\Omega <\infty.   \]
   The same claims extend to $\overline{\mathcal{A}^{1/2}(\Omega)}$. 
  \end{theorem}
  \begin{proof}
   Since $\Gamma$ is rectifiable, $\infty \notin \Gamma$. If $\infty \notin \Omega$, this is just  \cite[Theorem 10.3]{Duren_Hp_book}. If $\infty \in \Omega$, we can apply a M\"obius transformation $M$ so that $M(\Omega)$ is bounded
   and invoke conformal invariance of the norm and apply a change of variables. 
  \end{proof}
  
  In other words, the boundary values of elements of $\mathcal{A}^{1/2}(\Omega)$ and  $\overline{\mathcal{A}^{1/2}(\Omega)}$ are in $L^2(\Gamma,|dz|)$, that is, the $L^2$-space with respect to arc length.  
  
  By slightly adjusting the arguments in the proof above, or directly from the statement of the theorem using the polarization identity, we have 
  \begin{corollary} \label{co:inner_product_with_bvs}
  Let $\Omega$ be a simply connected domain in the sphere whose boundary $\Gamma$ is rectifiable.
   For $h_k\, dz^{1/2} \in \mathcal{A}^{1/2}(\Omega)$, $k=1,2$, we have   
   \[   \left( h_1 \,dz^{1/2},h_2 \,dz^{1/2} \right)_\Omega = \int_\Gamma h_1(z) \overline{h_2(z)} |dz|,   \]
   similarly for $\overline{h_k}\, d\bar{z}^{1/2} \in \overline{\mathcal{A}^{1/2}(\Omega)}$ we have 
   \[   \left( \overline{h_1} \, d\bar{z}^{1/2},\overline{h_2} \,d\bar{z}^{1/2} \right)_\Omega = \int_\Gamma \overline{h_1(z)} {h_2(z)} |dz|.   \]
  \end{corollary}
  \begin{remark} \label{re:bvs_fprime_exist}
   Since $\Gamma$ is rectifiable, if $f:\mathbb{D} \rightarrow \Omega$ is a conformal map then $f' \in H^1(\mathbb{D})$.
   In particular, $f'$ has non-tangential boundary values almost everywhere on $\partial \disk$ and 
   \[ \frac{d}{d\theta} f(e^{i\theta}) = i f'(e^{i\theta})  \]
   almost everywhere, where $f'(e^{i\theta})$ denote the non-tangential boundary values of $f'$.   
  \end{remark}
 
  Assume, for the moment, that the boundary $\partial \Omega$ of $\Omega$ is smooth.
  We then denote by $T: \partial \Omega \rightarrow S^{1} \subset \C$ the positively-oriented unit tangent vector.
  By \cite[Theorem 4.3]{Bell_book} every element $u$ of $L^2(\partial \Omega)$ has a unique decomposition
   \begin{equation} \label{eq:Bell_orthogonal_decomposition}
     u(\zeta)=h(\zeta) + \overline{H(\zeta)} \overline{T(\zeta)} 
   \end{equation}
   where $h$ and $H$ are non-tangential boundary values of elements of the Hardy space of $\Omega$.
   Furthermore this decomposition is orthogonal. As we observed above in Remark \ref{re:Smirnov_Hardy_confusion}, the Hardy space and Smirnov space agree in the case that the boundary is smooth, so that $h\, dz^{1/2} \in \mathcal{A}^{1/2}(\Omega)$ and $\overline{H} \, d\bar{z}^{1/2} \in \overline{\mathcal{A}^{1/2}(\Omega)}$.

\begin{remark}  In connection with Remark   \ref{re:Smirnov_Hardy_confusion}, we observe that although the decomposition (\ref{eq:Bell_orthogonal_decomposition}) is stated for Hardy spaces in S.~Bell's book \cite{Bell_book}, in fact it is more natural on Smirnov space. In particular, the appearance of the unit tangent vector $T$ is easily explained, and indeed the formula for the decomposition in the Smirnov space is more symmetric.    
 \end{remark}

\begin{remark}
    It would also make sense to write $h\sqrt{T} + \overline{H} \overline{\sqrt{T}}$ instead of $h + \overline{HT}$.
    This formulation makes the situation more apparently symmetric.
    However, we use the decomposition $h+ \overline{HT}$ because this is how it appears in the literature, moreover, the symmetry is not important for us at this point. (We establish some variation of the inherent symmetry later on in any case.)
\end{remark}
   
   One sees then that for smooth domains, one can identify $\mathcal{A}_h^{1/2}(\Omega)$ with $L^2(\Gamma,|dz|)$.
   The precise statement is the following theorem, which in fact only requires rectifiable boundary.
   Let 
   \[  W = \left\{ \left. h \right|_{\Gamma} \in L^{2}(\Gamma, |dz|) \, : h \, dz^{1/2} \in \mathcal{A}^{1/2}(\Omega) \right\} \]
   where $\left. h \right|_{\Gamma}$ denotes the non-tangential boundary values.  Then 
   \[ \overline{W} = \left\{ \left. \overline{H} \right|_{\Gamma} \in L^{2}(\Gamma, |dz|) :\overline{H}\, d\bar{z}^{1/2} \in \overline{\mathcal{A}^{1/2}(\Omega)} \right\}  \]
  is the set of complex conjugates of elements of $W$.

   \begin{theorem} \label{th:L2_decomposition}
    If $\Omega$ is a domain in the sphere with rectifiable boundary $\Gamma$, then we have the orthogonal decomposition
    \[  L^2(\Gamma,|dz|) = W \oplus \overline{T} \overline{W}.  \]
    Furthermore, the map 
    \begin{align*}
        \mathbf{b}_\Omega: \mathcal{A}_h^{1/2}(\Omega) & \rightarrow L^2(\Gamma,|dz|) \\
        h\, d{z}^{1/2} + \overline{H} d\bar{z}^{1/2} & \mapsto \left. h \right|_{\Gamma} + \left. \overline{H} \right|_{\Gamma} \overline{T}
    \end{align*}
    is an isometric isomorphism.
   \end{theorem}
   \begin{proof} 
Repeating the computation in Remark \ref{re:heuristic_inner_product} in the light of Theorem \ref{th:nontangential_bvs_exist} and Corollary \ref{co:inner_product_with_bvs}, shows that $\mathbf{b}_{\Omega}$ preserves the inner product.
Injectivity follows from orthogonality of $W \oplus \overline{T}\overline{W}$ together with Theorem \ref{th:nontangential_bvs_exist}.
Thus $\mathbf{b}_\Omega$ is an isometry.
It remains only to show that every element of $L^2(\Gamma,|dz|)$ is in $W \oplus \overline{T} \overline{W}$.  

Let $u \in L^2(\Gamma,|dz|)$. Then $u \circ f \sqrt{f'} \in L^2(\partial \disk,|dz|)$ (here we are using Remark \ref{re:bvs_fprime_exist}). It was already observed that for smooth domains $\Omega$, elements of $L^2(\partial \Omega,|dz|)$ have such a decomposition. In particular this holds for the disk, so 
\[   u \circ f \sqrt{f'} = h_1 + \overline{H_1} \overline{T_1} \]
for $h_1(z) dz^{1/2} \in \mathcal{A}^{1/2}(\disk), \overline{H_1(z)} d\bar{z}^{1/2} \in \overline{\mathcal{A}^{1/2}(\disk)}$, where $T_1(e^{i\theta})$ is the unit tangent vector on $\partial \disk$. 

By Remark  \ref{re:bvs_fprime_exist} the unit tangent vector on $\Gamma$ is 
\begin{equation*} 
 T(f(e^{i\theta})) = \frac{\frac{d f(e^{i\theta})}{d\theta}}{\left|\frac{d f(e^{i\theta})}{d\theta}\right|} = i e^{i\theta} \frac{\sqrt{f'(e^{i\theta})}}{\overline{\sqrt{f'(e^{i\theta})}}}
\end{equation*}
so since $T_1(e^{i\theta})=i e^{i\theta},$ we see that
    \begin{equation}     \label{eq:tangent_transformation}
     T \circ f = T_1  \frac{\sqrt{f'}}{\overline{\sqrt{f'}}}. 
   \end{equation}
 Thus we have
 \begin{equation*}
     u = \frac{h_1 \circ f^{-1}}{\sqrt{f' \circ f^{-1}}} + 
       \frac{H_1 \circ f^{-1}}{\sqrt{f' \circ f^{-1}}} \overline{T_1} \circ f^{-1}.
 \end{equation*}
 Setting 
 \[ h =  \frac{h_1 \circ f^{-1}}{\sqrt{f' \circ f^{-1}}}, \ \ \  H = \frac{H_1 \circ f^{-1}}{\overline{\sqrt{f' \circ f^{-1}}}} \]
 and applying (\ref{eq:tangent_transformation}) we obtain 
 \[  u  = h + \overline{H} \overline{T}.  \]
 The fact that $h dz^{1/2} \in \mathcal{A}^{1/2}(\Omega)$ and $\overline{H} d\bar{z}^{1/2} \in \overline{\mathcal{A}^{1/2}(\Omega)}$ follows from 
   Proposition \ref{pr:pull-back_is_isometry}. 
   \end{proof}
 \begin{remark} This theorem extends the decomposition \cref{eq:Bell_orthogonal_decomposition} given in \cite{Bell_book} to simply connected domains with rectifiable boundary. Although it is elementary, we were not able to locate this in the literature.
 \end{remark}
 The identification made in Theorem \ref{th:L2_decomposition} is motivated by the heuristic computation
   \begin{equation} \label{eq:heuristic_notation_change}
     h(\zeta) d\zeta^{1/2} + \overline{H(\zeta)}   d\overline{\zeta}^{1/2} =  \left( h(\zeta)  + \overline{H(\zeta)} \overline{T(\zeta)} \right) d\zeta^{1/2} 
   \end{equation}
    where we have used 
    \[  \overline{T} d\zeta^{1/2} = \overline{T} T^{1/2} |d\zeta|^{1/2} = 
    \overline{T}^{1/2} |d\zeta|^{1/2} = d\overline{\zeta}^{1/2}. \]
    Thus we obtain a function on $\Gamma$ from $u(\zeta) d\zeta^{1/2}$ by ``factoring out'' $d\zeta^{1/2}$. 
    
    This also suggests that we could write instead 
   \[  h(\zeta) d\zeta^{1/2} + \overline{H(\zeta)}   d\overline{\zeta}^{1/2} =  \left( h(\zeta) T(\zeta) + \overline{H(\zeta)}  \right) d\bar{\zeta}^{1/2}.
       \]
    Indeed one can show in exactly the same way that the orthogonal decomposition
    \begin{equation} \label{eq:alternate_L2_decomposition}
     L^2(\Gamma,|dz|) = T W \oplus \overline{W} 
    \end{equation} 
    holds. 
  
  Finally, we consider the following ``overfare''.
  Given $\mathbf{b}_\Omega$ as in Theorem \ref{th:L2_decomposition}, let 
  \[  \mathbf{b}^{-1}_\Omega:L^2(\Gamma,|dz|)  \rightarrow \mathcal{A}_h^{1/2}(\Omega)  \]
  denote its inverse. For a rectifiable Jordan curve $\Gamma$ in the sphere, let $\Omega_1$ and $\Omega_2$ be the two connected components of the complement.
  We then have that
  \[  \mathbf{b}^{-1}_{\Omega_2} \mathbf{b}_{\Omega_1}:\mathcal{A}_h^{1/2}(\Omega_1) \rightarrow \mathcal{A}^{1/2}_h(\Omega_2)           \]
  is a bounded map, which takes elements of $\mathcal{A}^{1/2}_h(\Omega_1)$ to elements of ${\mathcal{A}^{1/2}_h(\Omega_2)}$ with the same boundary values almost everywhere. 
  \begin{remark}
    An analogous ``overfare'' operator was defined by two of the authors in the case of Dirichlet spaces, under the assumption that the boundary is a quasicircle, see e.g.~\cite{Schippers_Staubach_Grunsky_expository}. In the Smirnov space setting the analysis is considerably more straightforward, at least for domains with rectifiable boundary.
  \end{remark}
  \begin{remark}
   If one uses the decomposition (\ref{eq:alternate_L2_decomposition}) to define $\mathbf{b}_{\Omega_1}$ and $\mathbf{b}^{-1}_{\Omega_2}$, then the resulting map does not change, so long as one makes a consistent choice on both sides.  To see this, assume that 
   \[  \mathbf{b}_{\Omega_1} \left( h(z) dz^{1/2} + \overline{H(z)} d\bar{z}^{1/2} \right) = \mathbf{b}_{\Omega_2} \left( h_*(z) dz^{1/2} + \overline{H_*(z)} d\bar{z}^{1/2} \right),  \]
   that is
   \begin{equation} \label{eq:overfare_welldefined_temp}
      h(z) + \overline{H(z)} \overline{T(z)} = h_*(z) + \overline{H_*(z)} \overline{T(z)}.  
   \end{equation}
   If we now define 
   \[  \tilde{\mathbf{b}}_{\Omega_1} \left( h(z) dz^{1/2} + \overline{H(z)} d\bar{z}^{1/2} \right) = h(z) T(z) + \overline{H(z)}   \]
   and similarly for $\tilde{\mathbf{b}}_{\Omega_2}$, multiplying both sides of (\ref{eq:overfare_welldefined_temp}) by $T$ we obtain that 
   \[  \tilde{\mathbf{b}}_{\Omega_1} \left( h(z) dz^{1/2} + \overline{H(z)} d\bar{z}^{1/2} \right) = \tilde{\mathbf{b}}_{\Omega_2} \left( h_*(z) dz^{1/2} + \overline{H_*(z)} d\bar{z}^{1/2} \right)  \]
  \end{remark}
  
  \begin{remark}
   The assumption of rectifiability can be weakened very slightly. It suffices to assume that $T(\Gamma)$ is rectifiable for some M\"obius transformation $T$.
   It is easily checked using the invariance of the norm and inner product that the results of this section can all be extended to this case. This might be of use in connection with the inversive geometry of \cite{Barrett_Bolt_Mobius,Barrett_Bolt_Laguerre}.   
  \end{remark} 
 \end{subsection}
 \end{section}
 \begin{section}{Faber operator and Faber series}\label{Faber Section}
 \begin{subsection}{Jump decomposition and Faber operators}
 
 We define the following Cauchy operator on $\overline{\mathcal{A}^{1/2}(\Omega)}$.
 Let $\Gamma$ be a rectifiable curve in the plane, and let $\Omega_1$ and $\Omega_2$ be the connected components of the complement in the sphere $\sphere$.  
 Given $\alpha \in {\mathcal{A}_h^{1/2}(\Omega_1)}$, define for $k=1,2$
 \begin{equation} \label{eq:Cauchy_first_definition}
   \left[ \mathbf{J}^{1/2}_{1,k}  \alpha  \right](z) =  \frac{1}{2 \pi i} \lim_{r \nearrow 1} \int_{\Gamma_{p,r}} 
    \alpha(w)\frac{dw^{1/2} dz^{1/2}}{w-z}, \ \ \ \ z\in \Omega_k.  
 \end{equation}
 Here we are following Barrett and Bolt \cite{Barrett_Bolt_Mobius} for the expression for the Cauchy operator and jump decomposition. 
 Their formalism of half-order differentials leads to an elegant (and computationally convenient) approach to the Cauchy integral in association with $L^2(\Gamma,|dz|)$ and $\mathcal{A}_h^{1/2}(\Omega_1)$ for rectifiable curves, as will be illustrated in the first part of this section. 
 A well-known result of G.~David (\cref{Davids}) says precisely when the jump decomposition holds.   
  \begin{remark}
     By Theorem \ref{th:L2_decomposition} we can also think of $\mathbf{J}_{1,k}^{1/2}$ as an operator on $L^2(\Gamma)$, where it is understood that the isometry $b_{\Omega_1}$ is used to identify the spaces $L^2(\Gamma)$ and $\mathcal{A}^{1/2}_h(\Omega_1)$.
  \end{remark}
    
 The integral is interpreted as follows. 
 Denoting $\alpha(z)=h(z) dz^{1/2} + \overline{H(z)} d\bar{z}^{1/2}$ we can write this as
 \begin{align} \label{eq:Cauchy_integral_breakdown}
   \left[ \mathbf{J}^{1/2}_{1,k} \alpha \right](z) & = \lim_{r \nearrow 1} \frac{1}{2 \pi i} 
 \int_{\Gamma_{p,r}} \frac{h(w)}{w-z}  dw \cdot dz^{1/2} + 
 \lim_{r \nearrow 1}\frac{1}{2 \pi i} \int_{\Gamma_{p,r}} 
    \frac{\overline{H(w)}}{w-z}  |dw| \cdot dz^{1/2} \ \ \ \ z\in \Omega_k \nonumber \\
    & =  \frac{1}{2 \pi i} 
 \int_{\Gamma} \frac{h(w)}{w-z}  dw \cdot dz^{1/2} + 
 \frac{1}{2 \pi i} \int_{\Gamma} 
    \frac{\overline{H(w)}}{w-z}  |dw| \cdot dz^{1/2} \ \ \ \ z\in \Omega_k  
 \end{align}
where in the second equality we are using Theorem \ref{th:nontangential_bvs_exist} to obtain the non-tangential boundary values of $h$ and $\overline{H}$ almost everywhere. 


 It is easily seen that the usual Cauchy integral on $L^2(\Gamma,|dz|)$ agrees with the integral (\ref{eq:Cauchy_integral_breakdown}) on $\mathcal{A}_h^{1/2}(\Omega_1)$ (still assuming $\Gamma$ is rectifiable). 
 For any $u \in L^2(\Gamma,|dz|)$ write $u= h + \overline{H} \overline{T}$ uniquely using Theorem \ref{th:L2_decomposition}. We then have 
 \begin{align*}
     \frac{1}{2\pi i} \int_{\Gamma} \frac{u(w)}{w - z} \, dw & =  \frac{1}{2 \pi i} \int_\Gamma \frac{h(w)}{w - z} dw + \int_{\Gamma} \frac{\overline{H(w)} \overline{T(w)}}{w-z} dw \\
     & = \frac{1}{2 \pi i} \int_\Gamma \frac{h(w)}{w - z} dw + \int_{\Gamma} \frac{\overline{H(w)}}{w-z} |dw|.
 \end{align*}
 Thus if we write $\alpha(z)=u(z) dz^{1/2}$ as in (\ref{eq:heuristic_notation_change}), the above integral agrees with (\ref{eq:Cauchy_integral_breakdown}). 
 
 In light of the above, we can remove the limit in equation (\ref{eq:Cauchy_first_definition}) to obtain
 \[  \left[ \mathbf{J}^{1/2}_{1,k}  \alpha  \right](z) =  \frac{1}{2 \pi i}  \int_{\Gamma} 
    \alpha(w)\frac{dw^{1/2} dz^{1/2}}{w-z} \ \ \ \ z\in \Omega_k.       \]
    and call it without reservation the Cauchy integral of $\alpha$. 
    
 \begin{proposition} \label{pr:Cauchy_integral_identity_on_holo}
  Let $\Gamma$ be a rectifiable curve in the plane, and let $\Omega_1$ and $\Omega_2$ be its complementary components in the Riemann sphere. Assume that $\infty \in \Omega_2$. For any $\alpha \in \mathcal{A}^{1/2}(\Omega_1)$ we have
  \[   \mathbf{J}^{1/2}_{1,2} \alpha    = 0  
  \ \ \ \ \mathrm{and} \ \ \ \      \mathbf{J}^{1/2}_{1,1} \alpha = \alpha.  \]
 \end{proposition}
 \begin{proof}
  This follows directly from (\ref{eq:Cauchy_integral_breakdown}). 
 \end{proof}

Using these jump operators, we now state a slight reformulation of David's result (Theorem \ref{Davids}).
By a result of M. Zinsmeister \cite{Zinsmeister}, Ahlfors-regular domains are Smirnov. Thus we obtain
 \begin{corollary} \label{co:Davids_theorem_halforder}
  Let $\Gamma$ be an Ahlfors-regular rectifiable Jordan curve of finite length and let $\Omega_1$ and $\Omega_2$ be the complementary components in the sphere.  Then
  \[  L^2(\Gamma,|dz|) = \mathcal{A}^{1/2}(\Omega_1) \oplus \mathcal{A}^{1/2}(\Omega_2). \]
 The decomposition is obtained from the bounded operators
 \[  \mathbf{J}^{1/2}_{1,k} :{\mathcal{A}_h^{1/2}(\Omega_1)} \rightarrow 
 \mathcal{A}^{1/2}(\Omega_k) \]
 for $k=1,2$.  
 \end{corollary}
 The decomposition above is the jump decomposition, and implicitly involves the non-tangential boundary values. More explicitly, any $H \in L^2(\Gamma,|dz|)$ which is boundary values of $\alpha \in \mathcal{A}_h^{1/2}(\Omega_1)$ satisfies
 \[   H = \mathbf{b}_{\Omega_1} \mathbf{J}_{1,1}^{1/2} \alpha - \mathbf{b}_{\Omega_2} \mathbf{J}_{1,2}^{1/2} \mathbf{b}_{\Omega_1} \alpha.      \]
 One can also write, for $\alpha \in \mathcal{A}_h^{1/2}(\Omega_1)$
 \begin{equation} \label{eq:overfared_jump}
    \alpha = \mathbf{J}_{1,1}^{1/2} \alpha - \mathbf{b}^{-1}_{\Omega_1} \mathbf{b}_{\Omega_2} \mathbf{J}_{1,2}^{1/2} \alpha.      
 \end{equation}

 We have the following immediate consequences of Corollary \ref{co:Davids_theorem_halforder}.
 \begin{corollary} \label{co:Nap_Yulm_analogue}
  Let $\Gamma$ be an Ahlfors-regular rectifiable Jordan curve of finite length and let $\Omega_1$ and $\Omega_2$ be the complementary components in the sphere. Then
  \[  \left. \mathbf{J}^{1/2}_{1,2} \right|_{\overline{\mathcal{A}^{1/2}(\Omega_1)}}:\overline{\mathcal{A}^{1/2}(\Omega_1)} \rightarrow   \mathcal{A}^{1/2}(\Omega_2) \]
  is an isomorphism.
 \end{corollary}
 \begin{proof}
  Assume that $\mathbf{J}^{1/2}_{1,2} \overline{H} d\bar{z}^{1/2} =0$. Set 
  \[   h(z) dz^{1/2} = -\mathbf{J}^{1/2}_{1,1}{\overline{H}(w)} d\bar{w}^{1/2}.  \]
  Then setting $\alpha = h(z) dz^{1/2} + \overline{H(z)} d\bar{z}^{1/2}$ it is easily checked that 
  \[  \left( \mathbf{J}^{1/2}_{1,1} \alpha, \mathbf{J}^{1/2}_{1,2} \alpha \right) = (0,0) \]
  so by Corollary \ref{co:Davids_theorem_halforder} $\alpha = 0$ so in particular $\overline{H(w)} d\bar{w}^{1/2}=0$. So $\mathbf{J}^{1/2}_{1,2}$ is injective on $\overline{\mathcal{A}^{1/2}(\Omega)}$. 
  
  Now let $g(z) dz^{1/2} \in \mathcal{A}^{1/2}(\Omega_2)$. 
  By Corollary \ref{co:Davids_theorem_halforder} there is an 
  \[  \alpha = h(z) dz^{1/2} + \overline{H(z)} d\bar{z}^{1/2} \]
  such that 
  \[  \left(\mathbf{J}^{1/2}_{1,1} \alpha, \mathbf{J}^{1/2}_{1,2} \alpha \right) = \left(0,g(z)dz^{1/2}\right).  \]
  By Proposition \ref{pr:Cauchy_integral_identity_on_holo} we see that 
  $\mathbf{J}^{1/2}_{1,2} \overline{H(w)} d\bar{w}^{1/2}= \mathbf{J}^{1/2}_{1,2} \alpha = g(z) dz^{1/2}$.  So $\mathbf{J}^{1/2}_{1,2}$ restricted to $\overline{\mathcal{A}_{1/2}(\Omega_1)}$ is surjective.
 \end{proof}
 This is an analogue in the Smirnov space setting of one direction a result of Napalkov and Yulmukhametov \cite{Nap_Yulm_Cauchy,Nap_Yulm}.\\
 
 Using the Cauchy operators we define the Faber operators as follows
 \begin{definition}\label{defn:faberiso}
  Let $\Gamma$ be an Ahlfors-regular rectifiable Jordan curve of finite length and let $\Omega_1$ and $\Omega_2$ be the complementary components in the sphere. Let $\f: \disk\to \Omega_1$ be a conformal map, and fix a choice of $\sqrt{\f'}$ to obtain an $\hat{\f} \in \Hdiffhat(\disk,\Omega_1)$. We define the \emph{Faber operator} as
 \[ \mathbf{I}_{\hat{\f}}^{1/2}=  -\mathbf{J}^{1/2}_{1,2} ({\hat{\f}}^{-1})^*: \overline{\mathcal{A}^{1/2}(\disk)} \rightarrow \mathcal{A}^{1/2}(\Omega_2). \]  
 \end{definition} 

 Thus we obtain the following theorem, see B. T. Bilalov and T. I. Najafov \cite{Bilalov_Najafov} where it appears with a different formulation.  
 \begin{corollary} \label{co:Faber_is_isomorphism} Let $\Gamma$ be an Ahlfors regular rectifiable Jordan curve. 
  The Faber operator is an isomorphism. 
 \end{corollary}
 \begin{proof}
  This follows directly from Proposition \ref{pr:pull-back_is_isometry} and Corollary \ref{co:Nap_Yulm_analogue}. 
 \end{proof}
 This is a counterpart in the Smirnov space setting of the result of A. \c{C}avu\c{s} \cite{Cavus} and Shen \cite{ShenFaber}, see Theorem \ref{th:Faber_operator_isomorphism} ahead.

 In fact by David's theorem, the converse holds under the much weaker assumption that the Faber operator is bounded, or equivalently, if $\mathbf{J}^{1/2}_{1,2}$ is bounded.  
 Assume that the boundary of the domain is rectifiable and $\mathbf{J}^{1/2}_{1,2}$ is bounded.  
 Given a rational function $R$ with no poles on $\Gamma$, it can be written as $R_{\Omega_1}+ R_{\Omega_2}$ where $R_{\Omega_k}$ is holomorphic on $\Omega_k$ for $k=1,2$. We have that $R$ is in $L^2(\Gamma)$.
 By David's theorem \cite[Theorem 3]{David} one has that if $\| R_{\Omega_2} \|_{L^2(\Gamma)} \lesssim  \| R \|_{L^2(\Gamma)}$ if and only if $\Gamma$ is Ahlfors regular. But $R_{\Omega_2} = \mathbf{J}^{1/2}_{1,2} R$, so since the Smirnov norm and $L^2$ norm on the boundary are comparable for rectifiable curves, this completes the proof.

 \end{subsection} 
 \begin{subsection}{Faber polynomials and Faber series}
 Let us very briefly recall some facts about the $p$-Faber series.  
 Let \(\Omega_1\) be a Jordan domain in \(\mathbb{C}\) with rectifiable boundary \(\Gamma\).
 Let \(\Psi: \mathbb{D}^{*}  \rightarrow \Omega_2\) be the Riemann mapping with \(\Psi^{\prime}(\infty)>0\), where $\Omega_2$ denotes the exterior of $\Omega_1$.
Let \(p'\) denote the H\"older-conjugate of \(p\geq 1\).
For \(k=0,1, \ldots\), and \(R>1\) define the polynomial of degree $k$
$$
\Phi_{p, k}(z)=\frac{1}{2 \pi i} \int_{|w|=R} \frac{w^{k}\left[\Psi^{\prime}(w)\right]^{1 / p'}}{\Psi(w)-z} d w, \quad z \in \Omega_1.
$$
These polynomials are referred to as $p$-Faber polynomials.\\
Now for any $g\in E^p(\Omega_1)$ one has
\begin{equation}\label{pfaber}
   g \sim \sum_{k=0}^{\infty} a_{k} \Phi_{p, k},
\end{equation}
where
\begin{equation}
    a_{k}=\frac{1}{2 \pi i} \int_{|w|=1} (g \circ \Psi(w)) \,\left[\Psi^{\prime}(w)\right]^{1 / p} w^{-k-1} d w,
\end{equation}
where  the series in \eqref{pfaber} is referred to as the $p$-Faber series of $g$.

For $p>1$, the $p$-Faber series for Jordan domains with rough boundary (i.e. boundary with corners) were studied by V. Kokila\v{s}vili \cite{Kokilasvili} and I. Ibragimov and D. Mamedhanov \cite{Ibramamed}. 
For Jordan domains with rectifiable boundary $p\geq 1$ the study of the $p$-Faber series was made by L--E. Andersson, who also gave a sufficient condition on the domain $\Omega$ for the bijectivity of the Faber operator 
\begin{equation}
\mathscr{F}_{p} \Phi(z)=\frac{1}{2 \pi i} \int_{|w|=1} \frac{\Phi(w) \,\left[\Psi^{\prime}(w)\right]^{1 / p}}{\Psi(w)-z} d w, \qquad z\in \Omega_1,
\end{equation}
when $p>1$ where $\Phi$ are boundary values of a function in $E^p(\disk)=H^p(\disk)$.  The main results necessary for our investigations, Theorem \ref{th:Faber_series_Smirnov} and Corollary \ref{co:Faber_is_isomorphism} above, are  reformulations of results in the Faber series literature; see \cite{Bilalov_Najafov} and references therein.   \\

 \c{C}avu\c{s} \cite{Cavus} and Shen \cite{ShenFaber} investigated the case of Faber series and operators, for $p=1$.
 One of the consequences of those investigations is the following result
  \begin{theorem} \label{th:Faber_operator_isomorphism}
  Let $\Gamma$ be a Jordan curve. Then the following are equivalent.
  \begin{enumerate}[font=\upshape]
      \item $\Gamma$ is a quasicircle.
      \item The Faber operator $\mathscr{F}_1 $ is a bounded isomorphism. 
    
  \end{enumerate}

 \end{theorem}
  
Now, returning to our study of half-order differentials, let $g_k(z) = \overline{z}^k d\bar{z}^{1/2}$.
  We define the $2$-Faber polynomial associated to the domain $\Omega_2$ via the conformal map $\f:\disk \rightarrow \Omega_1$ by 
  \[  \Phi_k = \mathbf{I}^{1/2}_{\hat{\f}} g_k   \]
  for $k \geq 0$, where the Faber operator $\mathbf{I}^{1/2}_{\hat{\f}}$ was introduced in Definition \ref{defn:faberiso}.  This is of the form 
  \begin{equation}\label{Phihat}
    \Phi_k = \hat{\Phi}_k(z) dz^{1/2}  
  \end{equation} 
  where $\hat{\Phi}_k(z)$ is a polynomial of degree $k+1$ in $1/z$. 
  
  In many of the sources on Faber series, including the ones cited above, the convention is that the Faber polynomials are defined on a bounded domain. Here we define them on the unbounded domain $\Omega_2$ in order to align with conventions in some of the literature on Teichm\"uller theory and in our previous papers \cite{Schippers_Staubach_Grunsky_expository}. 
  This change is inconsequential but for the convenience of the reader we note the change in the definition of the Faber operator, which with our convention would be defined by
  \begin{equation}
\mathscr{F}_{2} \Phi(z)= - \frac{1}{2 \pi i} \int_{|w|=1} \frac{\Phi(w) \,\left[\f^{\prime}(w)\right]^{1 / 2}}{\f(w)-z} d w, \qquad z\in \Omega_1,
\end{equation}
when $p>1$ where $\Phi$ are boundary values of a function in $E^2(\disk^*)$.  The change in sign is an artefact of the change in orientation induced by the change from the bounded to the unbounded side of the curve. 

The operators $\mathbf{I}_{\hat{\f}}^{1/2}$ and $\mathscr{F}_2$ are related by composition by $\mathbf{b}^{-1}_{\disk^*} \mathbf{b}_{\disk}$,
 which is just pull-back under the map $z \mapsto 1/\bar{z}$.  Namely, we have 
 \[ \mathscr{F}_2 \mathbf{b}^{-1}_{\disk^*} \mathbf{b}_{\disk} = \mathbf{I}^{1/2}_{\hat{\f}}.  \]
 
 \begin{remark}
 In the Dirichlet space setting, it has been found in earlier papers of the second two authors that this reformulation of the Faber operator on a space of anti-holomorphic functions leads to simpler functional analytic and function-theoretic identities. The formulations in this paper could be seen to confirm this, but we will not deal with this point in detail. 
 \end{remark}
 
   For any $h(z) dz^{1/2} \in \mathcal{A}^{1/2}(\Omega_2)$, let 
   \[ \overline{G(z)} d\bar{z}^{1/2} = {(\mathbf{I}^{1/2}_{\hat{\f}})}^{-1} h(w) dw^{1/2}.  \]
Then $\overline{G(z)} d \bar{z}^{1/2}$ is in the Smirnov space of the disk (which agrees with the Hardy space of the disk) and therefore the power series  of $\overline{G(z)}$ converges to $\overline{G(z)}$ in the Smirnov space. Denoting the power series by
   \[  \overline{G(z)} d\bar{z}^{1/2} = \sum_{k=0}^\infty a_k \bar{z}^k d\bar{z}^{1/2},  \] this leads to the definition of the Faber series of half-order differentials.
\begin{definition}
  We define the Faber series of $h(z) dz^{1/2}\in \mathcal{A}^{1/2}(\Omega_2)$ by 
  \[  \sum_{k=0}^\infty a_k \hat{\Phi}_k(z) dz^{1/2},  \]
  where
  $\hat{\Phi}_k(z)$ is given by \eqref{Phihat}.  
\end{definition}
   
   By applying the isomorphism ${\mathbf{I}^{1/2}_{\hat{\f}}}$ to the power series we obtain the following.
  \begin{theorem} \label{th:Faber_series_Smirnov}
   Let $\Gamma$ be a rectifiable Ahlfors-regular Jordan curve, and let $\Omega_1$ and $\Omega_2$ be the bounded and unbounded components of the complement respectively.  The Faber series of any element $h(z)\, dz^{1/2} \in \mathcal{A}^{1/2}(\Omega_2)$ converges in $\mathcal{A}^{1/2}(\Omega_2)$ to $h(z)\, dz^{1/2}$.  It is the unique series in $2$-Faber polynomials which does so.   
  \end{theorem} 
  \begin{remark}
      We note that this result is the counterpart of Theorem \ref{th:Faber_operator_isomorphism} above for the case of 2-Faber series.
  \end{remark}
 \end{subsection}
 \end{section}
 \begin{section}{Grunsky operator}
 \begin{subsection}{Szeg\H{o} and Garabedian kernels}
  In this section, we describe the Szeg\H{o} and Garabedian kernels in the half-order differential formalism of Barrett and Bolt/Hawley and Schiffer. 
  We will make use of an identity for the Garabedian kernel, which is analogous to an identity which appears in the work of Schiffer in the setting of Bergman spaces.  Although the Szeg\H{o} kernel is not necessary in the remaining results, we nevertheless included the reformulation, since it elucidates this point of view and fits quite naturally in the exposition.
 
 Let 
\begin{align*}
    \mathbf{P}_{\Omega}:\mathcal{A}_h^{1/2}(\Omega) &\rightarrow \mathcal{A}^{1/2}(\Omega), & &\text{ and }& \overline{\mathbf{P}_{\Omega}}:\mathcal{A}_h^{1/2}(\Omega) &\rightarrow \overline{\mathcal{A}^{1/2}(\Omega)}
\end{align*}
 denote the orthogonal projections.
 The decomposition is independent of the choice of point $p$; thus we use the notation $\mathbf{P}_{\Omega}$ rather than $\mathbf{P}(\Omega,p)$.
 Because the decomposition $\mathcal{A}_h^{1/2}(\Omega) = \mathcal{A}^{1/2}(\Omega) \overline{\oplus\mathcal{A}^{1/2}(\Omega)}$ is orthogonal we have $\mathbf{1} = \mathbf{P}_\Omega + \overline{\mathbf{P}_{\Omega}}$.
 
 This decomposition commutes with pull-back, and therefore if we have $\hat{g} \in \Hdiffhat(\Omega_{1},\Omega_{2})$ we immediately have that
 \begin{equation} \label{eq:projection_intertwines_pullback}
   \hat{g}^* \mathbf{P}_{\Omega_2} = \mathbf{P}_{\Omega_1} \hat{g}^*, \ \ \ \mathrm{and} \ \ \  \hat{g}^* \overline{\mathbf{P}_{\Omega_2}} = \overline{\mathbf{P}_{\Omega_1}} \hat{g}^*.   
 \end{equation}
 
 We will give integral expressions for the projections; these are the familiar Szeg\H{o} and Garabedian kernels. The Szeg\H{o} kernel for the disk is 
 \[  S_{\mathbb{D}}(z,\zeta) d\overline{\zeta}^{1/2} dz^{1/2} = \frac{1}{2\pi } \frac{d\overline{\zeta}^{1/2} dz^{1/2} }{1-\overline{\zeta} z}  \]
 and the Garabedian kernel is 
 \[   L_{\mathbb{D}}(\zeta,z) d {\zeta}^{1/2} dz^{1/2} = \frac{1}{2\pi } \frac{d {\zeta}^{1/2} dz^{1/2} }{\zeta -z}.  \]
 Let $\Omega$ be a simply-connected domain, and let $\hat{F} = (F,\sqrt{F'}) \in \Hdiffhat(\Omega, \disk)$.
 We define the Szeg\H{o} and Garabedian kernels of $\Omega$ by
 \begin{equation} \label{eq:Szego_definition}
   S_\Omega(z,\zeta) d \overline{\zeta}^{1/2} dz^{1/2} = (\hat{F} \times \hat{F})^* \left[ S_{\disk}(z,\zeta) d \overline{\zeta}^{1/2} dz^{1/2}\right] =\frac{1}{2\pi } \frac{\overline{\sqrt{F'(\zeta)}} \sqrt{F'(z)}}{1-\overline{F(\zeta)} F(z)}  d \overline{\zeta}^{1/2} dz^{1/2}  
 \end{equation}
 and 
 \begin{equation} \label{eq:Garabedian_definition}
   L_\Omega(\zeta,z) d\zeta^{1/2} dz^{1/2} = (\hat{F} \times \hat{F})^* \left[ L_{\disk}(\zeta,z) d {\zeta}^{1/2} dz^{1/2}\right] =\frac{1}{2\pi } \frac{\sqrt{{F'(\zeta)}} \sqrt{F'(z)}}{{F(\zeta)} - F(z)}  d {\zeta}^{1/2} dz^{1/2}.    
 \end{equation}
 These kernels do not depend on the choice of $(F,\sqrt{F'})$.
 Indeed, it is clear that the choice of sign of $\sqrt{F'}$ is immaterial.
 However, $F$ can still be replaced by $T \circ F$ for any disk automorphism $T$.
 For any M\"obius transformation $T$ it is easily verified that
 \[  \frac{\sqrt{T'(w)}\sqrt{T'(z)}}{T(w)-T(z)} = \frac{1}{w-z} \]
 and if $T$ is a disk automorphism, i.e. 
 \[  T(w)= e^{i\theta} \frac{w-a}{1-\bar{a}w} \]
 for some $a \in \disk$ then 
 \[  \frac{\overline{\sqrt{T'(w)}} \sqrt{T'(z)}}{1-\overline{T(w)} T(z)} = \frac{1}{1-\overline{w} z}  \]
 which establishes that the kernel functions are well-defined.

 With this definition, we immediately have that the Szeg\H{o} and Garabedian kernels are conformally invariant. That is, if $\hat{g} \in \Hdiffhat(\Omega_{1},\Omega_{2})$ then, if $z=g(w)$, $\eta=g(\zeta)$, we have
 \[ (\hat{g} \times \hat{g})^* \left(S_{\Omega_2}(z,\zeta) dz^{1/2} d\overline{\zeta}^{1/2} \right) =  S_{\Omega_1}(w,\eta) dw^{1/2} d\overline{\eta}^{1/2}   \]
 and
 \[ (\hat{g} \times \hat{g})^* \left(L_{\Omega_2}(z,\zeta) dz^{1/2} d{\zeta}^{1/2} \right) =  L_{\Omega_1}(w,\eta) dw^{1/2} d {\eta}^{1/2}.   \]
 Observe once again that the left hand sides are unchanged if one changes the choice of sign of $\sqrt{g'}$.
 
 We then have the following formulas for the projection operators.
 \begin{proposition} \label{pr:Szego_is_projection}
   Let $\Omega$ be a simply-connected domain which is conformally equivalent to the disk. We have for any $\alpha \in \mathcal{A}^{1/2}_h(\Omega)$ that,  
   \begin{equation}\label{eq:SzegoProjection}
        {\mathbf{P}_{\Omega}} \alpha = 
     \left( \alpha, \overline{S_\Omega(z,\zeta) d\overline{\zeta}^{1/2} dz^{1/2}} \right) =   \left( \lim_{r \nearrow 1} \int_{\Gamma_r}  S_\Omega(z,\zeta) d\overline{\zeta}^{1/2} \alpha \right) dz^{1/2}
    \end{equation}
  where $\Gamma_r=\{ z: g_p(z)=- \log{r} \}$ for some specific choice of $p$. 
 \end{proposition}
 \begin{proof}
  We prove that the formula holds in $\mathcal{A}^{1/2}(\Omega)$ and $\overline{\mathcal{A}^{1/2}(\Omega)}$ separately, and the result follows by orthogonality. If $\alpha= \overline{H(w)} d\overline{w}^{1/2}$ then since $\overline{H(\zeta)}$ and $S_\Omega(\zeta,z)$ are both anti-holomorphic in $\zeta$ we have 
  \[ \left( \alpha, \overline{S_\Omega(z,\zeta) d\overline{\zeta}^{1/2} dz^{1/2}} \right) =   \lim_{r \nearrow 1} \int_{\Gamma_r}  S_\Omega(z,\zeta)  \overline{H(w)} d\overline{\zeta} dz^{1/2} =0   \]
  and since ${\mathbf{P}_{\Omega}} \alpha = 0$ this proves the claim.
  
  Now assume that $\overline{\alpha} = h(\zeta) d\zeta^{1/2} \in \mathcal{A}^{1/2}(\Omega)$.
  Observe first that if $\Omega = \disk$ (or a sufficiently regular domain), \cref{eq:SzegoProjection} is a well-known fact in slightly different notation; see for example \cite[p.~23]{Bell_book}. 
  Let $F:\Omega \rightarrow \disk$ be a conformal map and $C_r= \{z: |z|=r \}$.  Then, denoting $\eta = F(\zeta)$ and $w=F(z)$, 
  \begin{align*}
    \left( \alpha, \overline{S_\Omega(z,\zeta) d\overline{\zeta}^{1/2} dz^{1/2}} \right) & = \lim_{r \nearrow 1} \int_{F^{-1}(C_r)} {S_\Omega(z,\zeta)  d\overline{\zeta}^{1/2} dz^{1/2} h(\zeta) d\zeta^{1/2}}   \\
    & = \lim_{r \nearrow 1} \int_{F^{-1}(C_r)}  \frac{\sqrt{F'(z)}\overline{\sqrt{F'(\zeta)}}}{1-\overline{F(\zeta)}F(z)}  d\overline{\zeta}^{1/2} dz^{1/2} h(\zeta) d\zeta^{1/2} \\
    & = \lim_{r \nearrow 1} \int_{C_r}  \frac{\sqrt{F'(z)}}{1-\overline{\eta}F(z)}  d\overline{\eta}^{1/2} dz^{1/2} h(F^{-1}(\eta)) \sqrt{(F^{-1})'(\eta)} d\eta^{1/2} \\
    & = F^* \left[ \lim_{r \nearrow 1} \int_{C_r} S_{\disk}(\eta,w) d\overline{\eta}^{1/2} dw^{1/2}  ((F^{-1})^* \alpha)(\eta)  \right] \\
    & = \left[ F^* \mathbf{P}_{\mathbb{D}} (F^{-1})^* \alpha \right](z).
  \end{align*}
  The claim now follows from \eqref{eq:projection_intertwines_pullback}.  
 \end{proof}
 \begin{proposition} \label{pr:Garabedian_projection}
  Let $\Omega$ be a simply-connected domain which is conformally equivalent to the disk. We have for any $\alpha \in \mathcal{A}_h^{1/2}(\Omega)$ that
  \[
  \left[\overline{P_{\Omega}} \alpha \right](\zeta) = 
     \left( \alpha,  {\frac{1}{i}L_\Omega(\zeta,z) d{\zeta}^{1/2} dz^{1/2}} \right) =   \lim_{r \nearrow 1} \int_{\Gamma_r} \overline{\frac{1}{i} L_\Omega(\zeta,z) d{\zeta}^{1/2} dz^{1/2}} \, \alpha    \]    
 \end{proposition}
 \begin{proof}
  By conformal invariance of the Garabedian kernel and the projection, a change of variable similar to the proof of Proposition \ref{pr:Szego_is_projection}, it suffices to prove this for any particular domain, for example the disk. 
  Again, we write $\alpha(\zeta)=h(\zeta) d\zeta^{1/2}+\overline{H(\zeta)}d\overline{\zeta}^{1/2}$.  
  
  First assume $\alpha(\zeta)=h(\zeta) d\zeta^{1/2}$.  It is immediately clear on the disk, using the Cauchy integral theorem, that 
  \[ \alpha(z)= H(z) dz^{1/2} = \frac{1}{i} \lim_{r \nearrow 1} \int_{|z|=r} L(\zeta,z) d\zeta^{1/2} dz^{1/2} H(\zeta) d\zeta^{1/2} = 
  \left(  \frac{1}{i} \lim_{r \nearrow 1} \int_{|z|=r} L(\zeta,z)  H(\zeta) d\zeta  \right) dz^{1/2}. \]
  Taking the complex conjugate proves the claim on the disk, which as noted above is sufficient. (Note that the general case for smooth domains is \cite[p25]{Bell_book}, if one sets $h=0$ in the decomposition of \cite[Theorem 4.3]{Bell_book}). 
  
  Now assume that $\alpha(\zeta)=h(\zeta) d\zeta^{1/2}$.  If the domain $\Omega$ is smoothly bounded, we have that
  \[ \frac{1}{i} \lim_{\epsilon \searrow 0} \int_{\Gamma_\epsilon}  L(\zeta,z) d\zeta^{1/2} dz^{1/2} \overline{h(\zeta)}d\overline{\zeta}^{1/2} =
  \left(\frac{1}{i} \lim_{\epsilon \searrow 0} \int_{\Gamma_\epsilon}  L(\zeta,z)  \overline{h(\zeta)} ds_\zeta \right)  dz^{1/2} =0 \]
  where the final equality is given in \cite[p25]{Bell_book}, if one sets $H=0$ in the decomposition of \cite[Theorem 4.3]{Bell_book}.
  Once again taking the complex conjugate proves the claim for smooth domains, which as observed above is sufficient.
 \end{proof}
 
 \begin{remark}
  The inner product with the Szeg\H{o} kernel can be written in the following way in terms of the contour integrals.  Letting $\alpha(\zeta)=h(\zeta) d\zeta^{1/2} + \overline{H(\zeta)}d\overline{\zeta}^{1/2}$, and $ds$ denotes infinitesimal arc length, we have
  \begin{align*}
     \left( \alpha, \overline{S_\Omega(z,\zeta) d\overline{\zeta}^{1/2} dz^{1/2}} \right) & = \lim_{\epsilon \searrow 1} \int_{\Gamma_\epsilon} S_\Omega(z,\zeta) d\overline{\zeta}^{1/2} dz^{1/2} \left(h(\zeta) d\zeta^{1/2} + \overline{H(\zeta)}d\overline{\zeta}^{1/2} \right) \\
     & = \lim_{\epsilon \searrow 1} \int_{\Gamma_\epsilon} S_\Omega(z,\zeta) d\overline{\zeta}^{1/2} dz^{1/2} h(\zeta) d\zeta^{1/2}  +\lim_{\epsilon \searrow 1} \int_{\Gamma_\epsilon}  S_\Omega(z,\zeta) d\overline{\zeta}^{1/2} dz^{1/2} \overline{H(\zeta)}d\overline{\zeta}^{1/2} \\
     & = \left( \lim_{\epsilon \searrow 1} \int_{\Gamma_\epsilon} S_\Omega(z,\zeta)  h(\zeta) ds_{\zeta} \right)  dz^{1/2} +   \left( \lim_{\epsilon \searrow 1} \int_{\Gamma_\epsilon} S_\Omega(z,\zeta) \overline{H(\zeta)} d\overline{\zeta} \right)  dz^{1/2}.
  \end{align*}
  Similarly, for the Garabedian kernel we have 
  \[ \left( \alpha,  {\frac{1}{i}L_\Omega(\zeta,z) d{\zeta}^{1/2} dz^{1/2}} \right) = \left( \lim_{\epsilon \searrow 1} \int_{\Gamma_\epsilon} \overline{\frac{1}{i} L(\zeta,z)} h(\zeta) ds_\zeta \right) d\overline{z}^{1/2} + \left( \lim_{\epsilon \searrow 1} \int_{\Gamma_\epsilon} \overline{\frac{1}{i} L(\zeta,z)} \overline{H(\zeta)} d\overline{\zeta} \right) d\overline{z}^{1/2}.  \]
 \end{remark}

 \end{subsection}
 \begin{subsection}{The Grunsky operator}

  We now define a generalization of the Grunsky operator to Smirnov spaces. 
  \begin{definition}\label{defn:Grunsky}
      Let $\Gamma$ be an Ahlfors-regular  rectifiable Jordan curve dividing the sphere into $\Omega_1$ and $\Omega_2$, and assume that $\infty \in \Omega_2$. Let $\f:\mathbb{D} \rightarrow \Omega_1$ be a conformal map.  We define the \emph{Grunsky operator} as the bounded operator given by 
  \begin{equation}\label{eq:HalfOrderGrunsky}
      \mathbf{Gr}_\f^{1/2} \defeq - {\hat{\f}}^* \, \mathbf{P}_{\Omega_1} \mathbf{b}^{-1}_{\Omega_1} \mathbf{b}_{\Omega_2} \, \mathbf{J}^{1/2}_{1,2} \, ({\hat{\f}}^{-1})^* \, : \, \overline{\mathcal{A}^{1/2}(\disk)} \rightarrow \mathcal{A}^{1/2}(\disk), 
  \end{equation}
  where we have chosen a branch of $\sqrt{\f'}$
  in order to obtain an element $\hat{\f} \in \Hdiffhat(\disk,\Omega_1)$.    
  \end{definition}
Since this is done consistently in the conjugation, $\mathbf{Gr}_\f^{1/2}$ is independent of this choice, justifying the notation. 
 \begin{proposition} \label{pr:Grunsky_J11_form}
  If $\Gamma$ is an Ahlfors-rectifiable Jordan curve, separating the sphere into $\Omega_1$ and $\Omega_2$, with $\infty \in \Omega_2$, then 
  \[   \mathbf{Gr}_\f^{1/2} = - \hat{\f}^* \mathbf{J}_{1,1}^{1/2} (\hat{\f}^{-1})^*.   \]
 \end{proposition}
 \begin{proof}
  Let $\overline{\alpha} \in \overline{\mathcal{A}^{1/2}(\disk)}$ and $\overline{\beta}=(\hat{\f}^{-1})^* \overline{\alpha}$. Using (\ref{eq:overfared_jump}) we obtain
  \begin{align*}
      \mathbf{Gr}_\f^{1/2} \overline{\alpha} & = -\hat{\f}^* \mathbf{P}_{\Omega_1} \mathbf{b}^{-1}_{\Omega_1} \mathbf{b}_{\Omega_2}\mathbf{J}_{1,2}^{1/2} \overline{\beta} \\
      & = - \hat{\f}^* \mathbf{P}_{\Omega_1} \left[ \mathbf{J}_{1,1}^{1/2} \overline{\beta} - \overline{\beta} \right] \\
      & = - \hat{\f}^* \mathbf{J}_{1,1}^{1/2} \overline{\beta}
  \end{align*}
  which proves the claim. 
 \end{proof}
  
 We also have the following integral expression for $\mathbf{Gr}_\f^{1/2}$ reminiscent of the Bergman-Schiffer integral expression for the Grunsky operator \cite{BergmanSchiffer,Schippers_Staubach_Grunsky_expository}.
 \begin{corollary} \label{co:Bergman_Schiffer_Grunsky}Let $\Gamma$ be a rectifiable Ahlfors-regular Jordan curve dividing the sphere into $\Omega_1$ and $\Omega_2$, and $\f:\mathbb{D} \rightarrow \Omega_1$ be a conformal map.  For any $\overline{\alpha} = \overline{H(z)} d\bar{z}^{1/2} \in \overline{\mathcal{A}^{1/2}(\disk)}$
 \begin{equation} \label{ferminonic Grunsky} [\mathbf{Gr}_\f^{1/2} \overline{\alpha}](\zeta) = 
  - \frac{1}{2 \pi i} \int_{\mathbb{S}^1,\eta} \left( \frac{\f'(\zeta)^{1/2} \f'(\eta)^{1/2}}{\f(\eta)-\f(\zeta)} - \frac{1}{\eta-\zeta}  \right)  d\eta^{1/2} d\zeta^{1/2} \overline{H(\eta)} d\bar{\eta}^{1/2}.   \end{equation}
 \end{corollary}
 \begin{proof}  Let $F=\f^{-1}$ and $\overline{\beta}= F^* \overline{\alpha}$. By (\ref{eq:Garabedian_definition}) and Proposition \ref{pr:Garabedian_projection} we have that 
 \[  \lim_{r \nearrow 1} \int_{\f(C_r),w} \frac{\sqrt{{F'(w)}} \sqrt{F'(z)}}{{F(w)} - F(z)}  d {w}^{1/2} dz^{1/2} \; \overline{\beta}(w) =0.     \]
 Thus 
 \begin{align*}
     \left[ \mathbf{J}_{1,1}^{1/2} \overline{\beta} \right](z) & = \frac{1}{2\pi i} \lim_{r \nearrow 1} \int_{\f(C_r),w} \left( \frac{1}{w-z} - \frac{\sqrt{F'(w)} \sqrt{F'(z)}}{F(w)-F(z)} \right) dw^{1/2} dz^{1/2} \; \overline{\beta}(w).     
 \end{align*}
 Now using Corollary \ref{pr:Grunsky_J11_form} and applying a change of variables, the result follows. 
 \end{proof}

 Of course, we could subtract any multiple of the Garabedian kernel and the result will still be true. The importance of the integral kernel in Corollary \ref{co:Bergman_Schiffer_Grunsky} is that it is non-singular. This trick of removing the singularity was applied in the setting of Bergman spaces by Schiffer \cite{Schiffer_first,BergmanSchiffer}.
 
  We also have the following result, which is an analogue of \cite[Theorems 6.5, 6.11]{Schippers_Staubach_Grunsky_expository}. 
 \begin{theorem}\label{thm:GraphOfGrunsky}  Let $\Gamma$ be a rectifiable Ahlfors-regular Jordan curve dividing the sphere into $\Omega_1$ and $\Omega_2$, and $\hat{\f} \in \Hdiffhat(\disk, \Omega_{1})$. Then 
   \begin{align*}
     \overline{\mathbf{P}_{\disk}}\, \hat{\f}^* \mathbf{b}^{-1}_{\Omega_1} \mathbf{b}_{\Omega_2} \mathbf{I}_{\hat{\f}}^{1/2} & = \mathbf{Id} \\
     \mathbf{P}_{\disk}\, \hat{\f}^* \mathbf{b}^{-1}_{\Omega_1} \mathbf{b}_{\Omega_2} \mathbf{I}_{\hat{\f}}^{1/2} & = \mathbf{Gr}_\f^{1/2}. 
   \end{align*}
   Thus, the graph of the Grunsky operator is the pull-back of the boundary values of $\mathcal{A}^{1/2}(\Omega_2)$. 
 \end{theorem}
 \begin{proof}
  The second claim follows almost immediately from the definition of the Grunsky operator, after observing that $\overline{\mathbf{P}_{\disk}}\,\hat{\f}^* = \hat{\f}^*\, \overline{\mathbf{P}_{\Omega_1}}$.
  The first claim follows from (\ref{eq:overfared_jump}), indeed
  Let $\overline{\alpha} \in \overline{\mc{A}^{1/2}(\disk)}$ and set $\overline{\beta} = (\hat{\f}^{-1})^* \overline{\alpha}$,
  \begin{align*}
  \overline{\mathbf{P}_{\disk}}\, \hat{\f}^* \mathbf{b}^{-1}_{\Omega_1} \mathbf{b}_{\Omega_2} \mathbf{I}_{\hat{\f}}^{1/2} \overline{\alpha} &=-
       \overline{\mathbf{P}_{\disk}} \hat{\f}^{*} \mathbf{b}_{\Omega_1}^{-1} \mathbf{b}_{\Omega_2} \mathbf{J}_{1,2}^{1/2} \overline{\beta} \\
      & = -\overline{\mathbf{P}_{\disk}} \hat{\f}^{*} \left[ \mathbf{J}_{1,1}^{1/2} \overline{\beta} - \overline{\beta} \right] \\
      & = \overline{\mathbf{P}_{\disk}} \hat{\f}^{*}  \overline{\beta} = \overline{\alpha}. \qedhere
  \end{align*}
 \end{proof}
 
 \end{subsection}
 \end{section}
 \begin{section}{The Weil-Petersson class and the Hilbert-Schmidt property of the Grunsky operator}\label{WP}

 In this section we investigate the relationship between the Hilbert-Schmidtness of the Grunsky operator and the Weil--Petersson class Teich""m\"uller space. More specifically, using the facts and the notations of Subsection \ref{Teichsection}, given a conformal map $\f:\disk \rightarrow \mathbb{C}$ onto a domain with rectifiable Ahlfors-regular boundary, we ask when the associated operator $\mathbf{Gr}_{\f}^{1/2}$ is an element of the Weil-Petersson class Teichm\"uller space.\\
 In the case of functions or 1-forms, it was shown independently by Takhtajan and Teo \cite{TNT}, and Shen \cite{Shen} that among conformal maps $\f^\mu$ onto domains bounded by quasicircles, the ``ordinary'' Grunsky operator $\mathbf{Gr}_{\f^{\mu}}$ is Hilbert-Schmidt if and only if is corresponding map $\f^{\mu}$ is in the Weil--Petersson class.

 One might hope for an equivalent statement concerning the Hilbert-Schmidt property of $\mathbf{Gr}_{\f}^{1/2}$ and the Weil-Petersson property of $\f$, but only one direction of the proof can be realized as a result for the operator $\mathbf{Gr}_{\f}^{1 / 2}$ itself.

 In this section, we assume throughout that $\Gamma$ is a rectifiable Ahlfors-regular Jordan curve dividing the sphere into $\Omega_1$ and $\Omega_2$, and $\f:\mathbb{D} \rightarrow \Omega_1$ is a conformal map.
 We consider the associated Grunsky operator $\mathbf{Gr}_{\f}^{1/2}: \overline{\mc{A}^{1/2}(\disk)} \rightarrow \mc{A}^{1/2}(\disk)$ as defined in \eqref{eq:HalfOrderGrunsky}.
 Let $\iota: \mathcal{A}^{1/2}(\disk) \to \mathcal{A}(\disk)$ be the inclusion operator.
 This operator is bounded, as a result of combination of the bounded inclusion \eqref{cont inclusion hardy-bergman} with the identification \eqref{Hardy_Ahalf correspondence}.
 We shall establish the following result.
\begin{theorem}\label{ WP-thm}  Let $\f:\disk \rightarrow \Omega$ be a conformal map. Assume that the boundary of $\Omega$ is Ahlfors-regular and rectifiable.
    The operator $\iota \mathbf{Gr}_{\f}^{1 / 2}$ is Hilbert-Schmidt if and only if $\f$ is in the Weil-Petersson class.
\end{theorem}
 This result is demonstrated by proving a couple of propositions and a lemma.\\ 
 
 We start by showing that the Hilbert-Schmidtness of $\iota \mathbf{Gr}_{\f}^{1 / 2}$ implies that $\f$ is in the WP--class, according to \textnormal{Definition \ref{defn:WP}}.
 In this connection, we also recall that an operator $A$ is Hilbert-Schmidt if and only if $AA^{\ast}$ is trace-class. Here $A^*$ denotes the adjoint of $A$. 

\begin{proposition}\label{trace yields WP}
    If $\iota\, \mathbf{Gr}_{\f}^{1 / 2}$ is Hilbert-Schmidt, then $\f$ is in the \textnormal{WP}--class.
\end{proposition}
\begin{proof}
The integral kernel of $ \iota\, \mathbf{Gr}_{\f}^{1 / 2}$ is given by
\begin{equation}\label{the original kernel}
    K(z, w)= \frac{1}{z-w}-\frac{\f^{\prime}(z)^{1 / 2} \f^{\prime}(w)^{1 / 2} }{\f(z)-\f(w)}
\end{equation}
and the kernel of $\Big(\iota \,\mathbf{G r}_{\f}^{1 / 2}\Big)^{\ast}$ is

\begin{equation}
    K^\ast(z, w)=\overline{\left(\frac{1}{w-z}-\frac{\f^{\prime}(z)^{1 / 2}\f^{\prime}(w)^{1 / 2} }{\f(w)-\f(z)}\right)}.
\end{equation}
From these it follows that the integral kernel of the operator $$\iota \mathbf{Gr}_{\f}^{1 / 2}\Big(\iota\, \mathbf{G r}_{\f}^{1 / 2}\Big)^{\ast}: \mathcal{A}(\disk)\to  \mathcal{A}(\disk) $$ is given by $M(z,w)$ with
\begin{align}\label{the kernel}
  &M(z,w)= \int_{\mathbb{S}^1} K(z,\zeta) K^{\ast}(\zeta, w) \, |d\zeta|  \\=\nonumber
  &\int_{\mathbb{S}^1} \left(\frac{1}{z-\zeta}- \frac{\f^{\prime}(z)^{1 / 2} \f^{\prime}(\zeta)^{1 / 2}}{\f(z)-\f(\zeta)}\right)\overline{\left(\frac{1}{w-\zeta}-\frac{\f^{\prime}(\zeta)^{1 / 2} f^{\prime}(w)^{1 / 2} }{\f(w)-\f(\zeta)}\right)} \, |d\zeta|.
  \end{align}
 This in turn yields that the trace of the operator $ \iota\, \mathbf{Gr}_{\f}^{1 / 2}\Big(\iota\, \mathbf{G r}_{\f}^{1 / 2}\Big)^{\ast}$ is given by
\begin{equation}
   \tau:= \int_{\mathbb{D}} M(z,z) \, dA_z= \int_{\mathbb{D}}\Big( \int_{\mathbb{S}^1} \left| \frac{\f^{\prime}(z)^{1/2} \f^{\prime}(\zeta)^{1/ 2}}{\f(\zeta)-\f(z)}-\frac{1}{\zeta-z}\right|^2 \, |d\zeta|\Big)\, dA_z.
    \end{equation}
    We assume now that $\iota \, \mathbf{Gr}_{\f}^{1/2}$ is Hilbert-Schmidt, and thus that $\tau < \infty$.
Our goal now is to show that this implies that $\f$ satisfies \eqref{Schwarzian WP condition},
which implies that $\f$ is in the {WP}-class.
Now suppose that
$$ \frac{1}{\zeta-z}-\frac{\f^{\prime}(z)^{1/2} \f^{\prime}(\zeta)^{1/ 2}}{\f(\zeta)-\f(z)}= \sum_{m,n=0}^{\infty} a_{m,n} z^n \zeta^m ,$$ where $a_{0,0}=0.$ Then one has
$$\int_{\mathbb{S}^1}  \Big(\frac{1}{\zeta-z}-\frac{\f^{\prime}(z)^{1/2} \f^{\prime}(\zeta)^{1/ 2}}{\f(\zeta)-\f(z)}\Big) \, |d\zeta|= \sum_{n=1}^{\infty} a_{0,n} z^n .$$

Note also that

$$\Big(\frac{1}{\zeta-z}-\frac{\f^{\prime}(z)^{1/2} \f^{\prime}(\zeta)^{1/ 2}}{\f(\zeta)-\f(z)}\Big)\Big|_{\zeta=0}= \sum_{n=1}^{\infty} a_{0,n} z^n . $$

Using a Taylor expansion, it can be shown that
\begin{equation}\label{schwarzian}
    \frac{1}{\zeta-z}-\frac{\f^{\prime}(z)^{1/2} \f^{\prime}(\zeta)^{1/ 2}}{\f(\zeta)-\f(z)}= C\, \mathcal{S}\f(z)(\zeta-z) +\mathrm{O}(|\zeta -z|^2),
\end{equation} for a constant $C$. This yields that

$$C\mathcal{S}\f(z)= \frac{1}{z}\Big(\frac{\f^{\prime}(z)^{1/2} \f^{\prime}(\zeta)^{1/ 2}}{\f(\zeta)-\f(z)}-\frac{1}{\zeta-z}\Big)\Big|_{\zeta=0} +\mathrm{O}(|z|)= -\sum_{n=1}^{\infty} a_{0,n} z^{n-1} + \mathrm{O}(|z|).$$

Therefore one has

 $$\int_{\mathbb{D}}(1-|z|^2)^2 |\mathcal{S}\f(z)|^2\, dA_z  \lesssim \sum_{n=1}^{\infty}  |a_{0,n}|^2 \int_{0}^{1} r^{2n-1} (1-r^2)^2 \, dr+ \mathrm{O}(1) \lesssim \sum_{n=1}^{\infty}  \frac{|a_{0,n}|^2}{n(n+1)(n+2)} +   \mathrm{O}(1)$$
 On the other hand a calculation shows that
 
$$ \tau= \frac{1}{2}\sum_{n,m=0}^{\infty}  \frac{|a_{n,m}|^2}{n+1} $$
and since
\begin{equation}
  \sum_{n=1}^{\infty}  \frac{|a_{0,n}|^2}{n(n+1)(n+2)} +   \mathrm{O}(1)  \lesssim \tau+ \mathrm{O}(1)
\end{equation}
we conclude that if $\tau<\infty$ then $\int_{\mathbb{D}}(1-|z|^2)^2 |\mathcal{S}\f(z)|^2\, dA_z <\infty.$
\end{proof}

\begin{corollary}\label{cor:HSImpliesWP}
    If $\mathbf{Gr}_{\f}^{1 / 2}$ is Hilbert-Schmidt, then $\f$ is in the Weil-Petersson class.
\end{corollary}
\begin{proof}
    Because $\iota$ is bounded, the assumption that $\mathbf{Gr}_{\f}^{1 / 2}$ is Hilbert-Schmidt implies that $\iota\, \mathbf{Gr}_{\f}^{1 / 2}$ is Hilbert-Schmidt (because of the well-known ``ideal-like'' property of Hilbert-Schmidt operators).
    By Prop.~\ref{trace yields WP} we then see that $\f$ is in the WP--class.
\end{proof}
We now proceed with the converse of Prop.~\ref{trace yields WP} i.e.~that the Weil-Petersson property yields Hilbert-Schmidtness.\\
To this end, we use the concepts that were introduced in Subsection \ref{Teichsection}.
 As was shown in \cite{TNT},
for $[\mu]\in T_{0}(1)$ with $\sup_{z\in \disk^*} \Big|\frac{\mu(z)}{\rho(z)}\Big|<\delta\in (0,1)$, $\delta$  sufficiently small, one can choose a representative $\mu\in L^2(\Del^*, \rho(z)\, dA_z)$ such that the path $[t\mu]$ connecting $0$ to $[\mu]$ in $T(1)$ lies in $T_{0}(1)$.  Let \(w_{t \mu}=\mathrm{g}_{t \mu}^{-1} \circ \mathrm{f}^{t \mu}\) be the corresponding conformal welding (according to \eqref{confweld}) and denote
by \(\left(K\right)_{t}(z, \zeta)\) the kernel \(K(z, \zeta)\) of $\iota \, \mathbf{G r}_{\mathrm{f^{t\mu}}}^{1 / 2}$ (given by \eqref{the original kernel}), associated with the conformal map \(\mathrm{f}^{t \mu}\). 
Then we obtain the following variational formula.
\begin{lemma}\label{varyK1}
Set $w_{t}=w_{t \mu},$ $\mathrm{f}_{t}=\mathrm{f}^{t \mu}$ and $\g_{t}=\mathrm{g}_{t \mu}.$ Then one has

\begin{equation}\label{viktig ekv}
  \left(\mu_{t} \circ \mathrm{g}_{t \mu}\right) \frac{\overline{\mathrm{g}_{t \mu}^{\prime}}}{\mathrm{g}_{t \mu}^{\prime}}=D_{t \mu} R_{(t \mu)^{-1}}(\mu).  
\end{equation}
and
\begin{equation}\label{eq:right_side_lemma}
\frac{d}{ds}\Big|_{s=0} (K)_{s+t}(\f_t^{-1}(z),
\f_t^{-1}(\zeta))\,\sqrt{(\f_t^{-1})'(z)}\,\sqrt{(\f_t^{-1})'(\zeta)}=\frac{-1}{2\pi}\int\limits_{ \Omega_t^*}
\frac{\mu_t(u)(\zeta-z)}{(u-z)^2(u-\zeta)^2}\, dA_u,
\end{equation} 
where \(\Omega_{t}^{*}=\mathrm{f}^{t \mu}\left(\mathbb{D}^{*}\right)=\mathrm{g}_{t \mu}\left(\mathbb{D}^{*}\right)\). 
\end{lemma}
Note that the integral in  \eqref{eq:right_side_lemma} is not a principal value, because $\zeta$ and $z$ are in $\Omega_t$.  
\begin{proof}  Since the equality \eqref{viktig ekv} was proven in Theorem 2.6 in \cite{TNT}, it only remains to show \eqref{eq:right_side_lemma}. Setting \(v_{s}=\mathrm{f}_{s+t} \circ \mathrm{f}_{t}^{-1}\),
the variational formula for quasiconformal mappings (see \cite{Ahlfors}) yields that
\begin{align}\label{variation}
\left.\frac{d}{ds}\right\vert_{s=0} v_s(z)= \dot{v}_0(z)=
-\frac{1}{\pi}\int\limits_{\Omega_{t}^*}
\frac{\mu_t(u)z(z-1)}{(u-z)u(u-1)}dA_u+ p(z),
\end{align}
where $p(z)$ is a quadratic polynomial. On the other hand
\begin{align}\label{goer1}
&(K)_{s+t}\left(\f_t^{-1}(z),
\f_t^{-1}(\zeta)\right)\sqrt{\left(\f_t^{-1}\right)'(z)\left(\f_t^{-1}\right)'(\zeta)}\, \\
=& \nonumber\sqrt{\left(\f_t^{-1}\right)'(z)\left(\f_t^{-1}\right)'(\zeta)}\, \left(\frac{1}{\f_t^{-1}(z)-\f_t^{-1}(\zeta)}-\frac{\sqrt{\f'_{s+t}(\f_t^{-1}(z))}\sqrt{\f'_{s+t}(\f_t^{-1}(\zeta))}
}{\f_{s+t}(\f_t^{-1}(z))-\f_{s+t}(\f_t^{-1}(\zeta))}\right)\\= &\nonumber
\sqrt{\left(\f_t^{-1}\right)'(z)\left(\f_t^{-1}\right)'(\zeta)}\left(\frac{1}{\f_t^{-1}(z)-\f_t^{-1}(\zeta)}-\frac{\sqrt{\f'_{s+t}(\f_t^{-1}(z))}\sqrt{\f'_{s+t}(\f_t^{-1}(\zeta))}
}{\f_{s}(
z)-\f_{s}(\zeta)}\right)\\= &\nonumber
\left(\frac{\sqrt{\left(\f_t^{-1}\right)'(z)\left(\f_t^{-1}\right)'(\zeta)}}{\f_t^{-1}(z)-\f_t^{-1}(\zeta)}-\frac{\sqrt{v'_{s}(z)\,v'_{s}(\zeta)}}{v_{s}(
z)-v_{s}(\zeta)}\right),
\end{align}
and
\begin{align*}
\left.\frac{d}{ds}\right|_{s=0}
\frac{\sqrt{v_s'(z)v_s'(\zeta)}}{v_s(z)-v_s(\zeta)}=\frac{(z-\zeta)(\dot{v}_0'(z)+\dot{v}_0'(\zeta))/2-(\dot{v}_0(z)-\dot{v}_0(\zeta))}{(z-\zeta)^2}.
\end{align*}
This yields that
\begin{align}\label{goer2}
\begin{split}
\frac{d}{d s}\Big|_{s=0} (K)_{s+t}\left(\f_t^{-1}(z), \f_t^{-1}(\zeta)\right)\sqrt{\left(\f_t^{-1}\right)'(z)}\,\sqrt{\left(\f_t^{-1}\right)'(\zeta)} \\
=-\frac{(z-\zeta)(\dot{v}_0'(z)+\dot{v}_0'(\zeta))/2-(\dot{v}_0(z)-\dot{v}_0(\zeta))}{(z-\zeta)^2}.
\end{split}
\end{align}

Now note that
\begin{align}\label{goer3}
\begin{split}
&\frac{1}{2(z-\zeta)}(\dot{v}_0'(z)+\dot{v}_0'(\zeta))\\=&-\frac{1}{2\pi}\int\limits_{\Omega_{t}^*}
\frac{\mu_t(u)}{u(u-1)(z-\zeta)}\Big(\frac{2uz-u-z^2}{(u-z)^2}+\frac{2u\zeta-u-\zeta^2}{(u-\zeta)^2}\Big)\, dA_u+ \frac{p'(z)+p'(\zeta)}{2(z-\zeta)}
\end{split}
\end{align}
and
\begin{equation}\label{goer4}
   -\frac{\dot{v}_0(z)-\dot{v}_0(\zeta)}{(z-\zeta)^2}= \frac{1}{\pi}\int\limits_{\Omega_{t}^*}
\frac{\mu_t(u)}{u(u-1)(z-\zeta)^2}\Big(\frac{z^2-z}{u-z}-\frac{\zeta^2-\zeta}{u-\zeta}\Big)\, dA_u - \frac{p(z)-p(\zeta)}{(z-\zeta)^2}.
\end{equation}
We also observe that
\begin{align}\label{goer5}
\begin{split}
    \frac{-1}{2(z-\zeta)u(u-1)}\Big(\frac{2uz-u-z^2}{(u-z)^2}+\frac{2u\zeta-u-\zeta^2}{(u-\zeta)^2}\Big)+\frac{1}{(z-\zeta)^2 u(u-1)}\Big(\frac{z^2-z}{u-z}-\frac{\zeta^2-\zeta}{u-\zeta}\Big)\\=\frac{\zeta-z}{2(u-z)^2(u-\zeta)^2}.
    \end{split}
\end{align}
This yields that
\begin{align}\label{skunt}
\begin{split}
    \frac{(z-\zeta)(\dot{v}_0'(z)+\dot{v}_0'(\zeta))/2-(\dot{v}_0(z)-\dot{v}_0(w))}{(z-\zeta)^2}\\= \frac{1}{2\pi}\int\limits_{ \Omega_t^*}
\frac{\mu_t(u)(\zeta-z)}{(u-z)^2(u-\zeta)^2}\, dA_u +  \frac{p'(z)+p'(\zeta)}{2(z-\zeta)}-\frac{p(z)-p(\zeta)}{(z-\zeta)^2}.
    \end{split}
\end{align}
However, since for any quadratic polynomial $p(z)= az^2 +bz+c$ one has that

\begin{align}
\begin{split}
    \frac{p'(z)+p'(\zeta)}{2(z-\zeta)}-\frac{p(z)-p(\zeta)}{(z-\zeta)^2}=\frac{(z-\zeta)(p'(z)+p'(\zeta))-2(p(z)-p(\zeta))}{2(z-\zeta)^2}\\= \frac{(z-\zeta)(2az+2a\zeta+2b)-2(az^2+bz-a\zeta^2 -b\zeta)}{2(z-\zeta)^2}=0,
    \end{split}
\end{align}

\eqref{goer2} and \eqref{skunt}  yield the desired result.\\

\end{proof}

Having the lemma at our disposal, we can prove that the operator $\iota \mathbf{Gr}_{\f^{\mu}}^{1 / 2}$ is Hilbert-Schmidt, if $\f^\mu$ is in the Weil-Petersson class.
\begin{proposition}\label{traceclass}
If  $\f^{\mu}$ and $\g_{\mu}$ correspond to a point $[\mu]$ in the \emph{WP}--class Teichm\"uller space $T_{0}(1)$, then the operator $\iota\, \mathbf{Gr}_{\f^{\mu}}^{1 / 2}\Big(\iota\, \mathbf{G r}_{\f^{\mu}}^{1 / 2}\Big)^{\ast}: \mathcal{A} (\disk) \to \mathcal{A}(\disk),$ is trace-class.
\end{proposition}
\begin{proof}
Recalling that \(\left(K\right)_{1}(z, \zeta):= (K)_{t=1}(z, \zeta)\) is the kernel \(K(z, \zeta)\) of $\iota \, \mathbf{G r}_{\mathrm{f}}^{1 / 2}$ (given by \eqref{the original kernel}), and reasoning as in the Proof of Proposition \ref{trace yields WP}, it is enough to show that for $\f^\mu=\f$ one has that
\begin{equation}\label{eq1}
   \int_{\mathbb{D}}\int_{
\mathbb{S}^1} \bigl|(K)_1(\zeta,z)\bigr|^2 \, |d\zeta|\, dA_z=\int_{\mathbb{D}}\Big( \int_{\mathbb{S}^1} \Big| \frac{\f^{\prime}(z)^{1/2} \f^{\prime}(\zeta)^{1/ 2}}{\f(\zeta)-\f(z)}-\frac{1}{\zeta-z}\Big|^2 \, |d\zeta|\Big)\, dA_z
    <\infty.
\end{equation}

To this end, by the fundamental theorem of calculus, Jensen's inequality and Fubini's theorem we have that
\begin{align}\label{spoor}
&\int_{\mathbb{D}}\int_{
\mathbb{S}^1} \bigl|(K)_1(\zeta,z)\bigr|^2 \, |d\zeta|\, dA_z\\ \nonumber
=&\int\limits_{\Del}\int\limits_{\mathbb{S}^1}\left|
 \int_0^1 \frac{d}{dt}(K)_t(\zeta,z)dt\right|^2\, |d\zeta|\, dA_z \\
\leq \nonumber & \int_0^1 \int\limits_{\Del}\int\limits_{\mathbb{S}^1}\left|
 \frac{d}{dt}(K)_t(\zeta,z)\right|^2\, |d\zeta|\, dA_z dt \\
= \nonumber &\int_0^1 \int\limits_{\Del}\int\limits_{\mathbb{S}^1}\left|
\frac{d}{ds}\Bigr\vert_{s=0}(K)_{t+s}(\zeta,z)\right|^2\, |d\zeta|\, dA_z dt\\ 
= \nonumber &\int_{0}^{1}I(t)dt.
\end{align}
Therefore, our goal is to show that $\int_{0}^{1}I(t)dt<\infty.$
A change of variables $\zeta\mapsto \f^{-1}_{t}(\zeta), z\mapsto \f^{-1}_{t}(z)$ in the inner integral $I(t),$ and Lemma \ref{varyK1}, yield that 
\begin{align*}
I(t)=
 & \int\limits_{\Omega_t}\int\limits_{\partial\Omega_t} \left|
\left.\frac{d}{ds}\right|_{s=0} (K)_{t+s}\left(\f_t^{-1}(\zeta),
\f_t^{-1}(z)\right)\sqrt{\left(\f_t^{-1}\right)'(\zeta)\left(\f_t^{-1}\right)'(z)}\right|^2 |\left(\f_t^{-1}\right)'(z)|
\, |d\zeta|\, dA_z\\
=&\frac{1}{2\pi}
\int\limits_{\Omega_t}\int\limits_{\partial\Omega_t} \left|\,
\int\limits_{ \Omega_t^*}
\frac{\mu_t(u) }{(u-\zeta)^2(u-z)^2}\, dA_u\,\right|^2\, |z-\zeta|^2 |\left(\f_t^{-1}\right)'(z)|  |d\zeta|\, dA_z.
\end{align*}
Now assume for the moment that $\sup_{z\in \Omega_t}|z-\zeta| |\left(\f_t^{-1}\right)'(z)|<\infty$ for $\zeta\in \partial\Omega_t$, which will be motivated at the end of the proof.
Then one has that
\begin{align*}
I(t)\lesssim & 
\int\limits_{\Omega_t}\int\limits_{\partial\Omega_t} \left|\,
\int\limits_{ \Omega_t^*}
\frac{\mu_t(u) }{(u-\zeta)^2(u-z)^2}\, dA_u\,\right|^2 \, |z-\zeta|\, |d\zeta|\, dA_z \\&\lesssim \int\limits_{\mathbb{C}}\int\limits_{\partial\Omega_t} \left|\,
\int\limits_{\mathbb{C}}
\frac{\chi_{\Omega^{*}_t}(u)\mu_t(u) }{(u-\zeta)^2(u-z)^2}\, dA_u\,\right|^2 \, |z-\zeta|\, |d\zeta|\, dA_z\\&\lesssim  \int_{\partial\Omega _t}\Big(\int_{\mathbb{C}}\left|\
\mathscr{B}F(z)\,\right|^2 \, |z-\zeta|\, dA_z \Big)\, |d\zeta|
\end{align*}
where  $F(z)= \frac{\chi_{\Omega^{*}_t}(z)\mu_t(z) }{(z-\zeta)^2},$ $\chi_{\Omega^{*}_t}$ is the characteristic function of $\Omega^{*}_t ,$
and 
$\mathscr{B} F(z)$ is the Beurling transform of $F$ given by \eqref{Beurling}. By Theorem \ref{CF theorem}, the Beurling transform (being a Calder\'on-Zygmund operator in $\R^2$), satisfies the weighted norm inequality
\begin{equation}\label{wnineq}
\int_{\mathbb{C}} |\mathscr {B} F(z)|^2 |z-\zeta|\, dA_z \lesssim  \int_{\mathbb{C}} | F(z)|^2 |z-\zeta|\, dA_z,
\end{equation}
because the Beurling transform commutes with translations, $|z|$ is a Muckenhoupt $A_2$-weight in $\mathbb{C},$ (see the facts in Subsection \ref{Muckweights} regarding power-weights), and $\mathscr {B}$ is bounded on the weighted space
$L^2_{|z|}(\mathbb{C}, dA_z).$ 
Therefore using \eqref{wnineq} we obtain
\begin{align*}
I(t)\leq &\frac{1}{2\pi}\int\limits_{\partial\Omega_{t}}\int\limits_{\Omega_t^*} \frac{ |
\mu_t(z)|^2}{|z-\zeta|^4} |z-\zeta| \, dA_z \, |d\zeta| =\int\limits_{\Omega_t^*} |
\mu_t(z)|^2 \Big(\int\limits_{\partial\Omega_{t}} \frac{|d\zeta|}{|z-\zeta|^{3}}\Big) \, dA_z \,  . 
\end{align*}

Note that if $\eta(z)= d(\partial \Omega_t ,z)$ (here $d(z,w):=|z-w|$), then one has 
$\eta(z)\sim (\rho_2)_{t}^{-1/2} (z)$ for $z\in \Omega^{*}_t$ (see \cite{Lehto}, \cite{Nagbook} for a proof) where $(\rho_2)_t (z) := (\rho\circ \g_t^{-1})(z) |(\g_t^{-1})'(z)|^2$ is the hyperbolic metric density on $\Omega^{*}_{t}$. Therefore 
\begin{equation} \label{inequality}
\int\limits_{\partial \Omega_{t}}\frac{|d \zeta|}{|z-\zeta|^{3}}\leq \int\limits_{|z-\zeta|=\eta(z)}\frac{|d \zeta|}{|z-\zeta|^{3}}= \frac{2\pi}{\eta^{2}(z)}\lesssim  (\rho_{2})_{t}(z),\,\,z\in\Omega^{*}_{t},
\end{equation}
This and \eqref{norm of mu} in turn yield that
\begin{align*}
I(t)\lesssim \int\limits_{\Omega_t^*}
|\mu_t(z)|^2 (\rho_2)_t(z)\, dA_z =
\int\limits_{\Del^*}|\tilde{\mu}_t(z)|^2 \rho(z)\, dA_z =
\Vert
\tilde{\mu}_t\Vert_2^2,
\end{align*}
where $\tilde{\mu}_t = D_{t\mu}R_{(t\mu)^{-1}}(\mu)$ and $R$ is the right translation in \eqref{defn:righttrans}, and we have also used \eqref{viktig ekv}.
Now it follows from Remark 2.8 in Chapter 1 of \cite{TNT} that 
\begin{align*}
\Vert \tilde{\mu}_t\Vert_2 \lesssim \Vert \mu \Vert_2 ,
\end{align*}
for all $0\leq t\leq 1$. This together with \eqref{spoor} finally yields that $\int_{0}^{1}I(t)dt\lesssim \Vert \mu\Vert^2_2$, which by Definition \ref{defn:WP} would conclude the proof.\\
Therefore it only remains to confirm the claim that $\sup_{z\in \Omega_t}|z-\zeta| |\left(\f_t^{-1}\right)'(z)|<\infty$ for $\zeta\in \partial\Omega_t$. Note that for $z\in \Omega_t$ one has (see e.g. \cite{Lehto} and \cite{Nagbook} for a proof) 
$$d(\partial \Omega_t,z)\sim (\rho_1)_{t}^{-1/2} (z):= (\rho\circ \f_t^{-1})^{-1/2}(z) |(\f_t^{-1})'(z)|^{-1}.$$This yields that for $\zeta\in \partial \Omega_t$
$$\sup_{z\in \Omega_t}|z-\zeta| |\left(\f_t^{-1}\right)'(z)|\lesssim \sup_{z\in \Omega_t} (\rho\circ \f_t^{-1})^{-1/2}(z)= \sup_{w\in \mathbb{D}}\frac{1}{2}(1-|w|^2)<\infty,$$ 
which shows the claim.
\end{proof}

Thm.~\ref{ WP-thm} gives necessary and sufficient conditions on $\f$, such that $\iota \mathbf{Gr}_{\f}^{1/2}$ is Hilbert-Schmidt, under the assumption that $\f$ is a conformal map from the disk to a domain with Ahlfors-regular and rectifiable boundary.
However, in the case of 1-forms mentioned at the beginning of this section, it is the operator $\mathbf{Gr}^{1/2}_{\f}$ that is studied.
Indeed, from the point of view of the applications, this is the operator that is of interest.
We are thus led to the following problem.
\begin{problem}\label{prob:characterization}
    Determine necessary and sufficient conditions on $\Omega$ and $\f: \disk \rightarrow \Omega$ which guarantee that the operator $\mathbf{Gr}^{1/2}_{\f}$ is Hilbert-Schmidt.
\end{problem}

Of course, for the operator $\mathbf{Gr}_{\f}^{1/2}$ to be Hilbert-Schmidt, it must first be bounded.
We have found that a sufficient condition for this, is that the boundary of $\Omega$ is Ahlfors-regular and rectifiable.
Under this assumption, for $\mathbf{Gr}_{\f}^{1/2}$ to be Hilbert-Schmidt, it is then necessary that $\f$ is Weil-Petersson (Cor.~\ref{cor:HSImpliesWP}).
A reasonable next step towards a resolution of Problem \ref{prob:characterization} would be to determine if this condition is also sufficient.

 \end{section}

 \appendix

 \begin{section}{Differential geometry of half-order differentials}\label{sec:HalfOrderGeo}
    In this section, we take a differential geometric approach to defining half-order differentials, essentially following Barrett and Bolt \cite{Barrett_Bolt_Mobius}.
    We shall use the same notation as in \cref{sec:HalfOrderDifferentials}, we shall see later that this is justified.

    Given a simply-connected domain $\Omega$, a (holomorphic) half-order differential should be an object that can be squared to give a differential of rank $(1,0)$ on $\Omega$. (Or more generally, one should be able to pair two different half-order differentials to obtain a $(1,0)$-differential.)

    We consider the two charts on the Riemann sphere $\phi_{0}:\sphere \setminus \{\infty \} \rightarrow \C, z \mapsto z$ and $\phi_{\infty}:\sphere \setminus \{ 0 \} \rightarrow \C, z \mapsto 1/z$.
    Let $E$ be the complex line bundle defined by the transition function $\sphere \supset \C^{\times} \xrightarrow{g_{0\infty}} \C^{\times} = \GL(\C), z \mapsto i/z$.
    The bundle $E \rightarrow \sphere$ is holomorphic, because $z \mapsto i/z$ is.
    We recall that the holomorphic cotangent bundle $T^{*}_{1,0}\sphere$ has transition function $\C^{\times} \rightarrow \C^{\times}, z \mapsto -1/z^{2}$.
    This implies that there is an isomorphism of holomorphic vector bundles $E \otimes E \cong T_{1,0}^{*}\sphere$.

    Let $\Omega \subseteq \sphere$ be a (not necessarily proper) subset.
    A \emph{holomorphic half-order differential} (or $\frac{1}{2}$-differential) on $\Omega$ is a holomorphic section of $E|_{\Omega}$.
    We denote the set of holomorphic half-order differentials on $\Omega$ by $\Omega^{\frac{1}{2},0}(\Omega)$.
    Similarly, one obtains an anti-holomorphic line bundle $\overline{E}$, together with an isomorphism $\overline{E} \otimes \overline{E} \rightarrow T^{*}_{0,1}\sphere$.
    Anti-holomorphic sections of $\overline{E}$ are called anti-holomorphic half-order differentials, they constitute $\Omega^{0,\frac{1}{2}}(\Omega)$.
    
    Using the identification $E \otimes E \cong T^{*}_{1,0}\sphere$ one may pair half-order differentials to obtain ordinary differentials, i.e.~we have bilinear maps
    \begin{align*}
        \Omega^{\frac{1}{2},0}(\Omega) \times \Omega^{\frac{1}{2},0}(\Omega) &\rightarrow \Omega^{1,0}(\Omega), & \Omega^{0,\frac{1}{2}}(\Omega) \times \Omega^{0,\frac{1}{2}}(\Omega) &\rightarrow \Omega^{0,1}(\Omega),
    \end{align*}
    which justifies the terminology.
    We denote the above maps simply by juxtaposition.

    If $\Omega \subset \sphere$ is an open such that $\infty \notin \Omega$, we write $dz \in \Omega^{1,0}(\Omega)$ for the section associated to the ``identity coordinate'' $\Omega \rightarrow \C, z \mapsto z$.
    Any element of $\Omega^{1,0}(\Omega)$ can be written as $h dz$, where $h: \Omega \rightarrow \C$ is a holomorphic function.
    By construction, the bundle $E|_{\Omega}$ comes equipped with a trivialization.
    We denote by $\sqrt{dz} \in \Omega^{\frac{1}{2},0}(\Omega)$ the corresponding holomorphic section.
    Any element of $\Omega^{\frac{1}{2},0}(\Omega)$ can be written as $h\sqrt{dz}$, where $h: \Omega \rightarrow \C$ is a holomorphic function.
    Observe that this notation is internally consistent, in the sense that 
    \begin{equation}\label{eq:pairing}
        h_{1}\sqrt{dz} \otimes h_{2}\sqrt{dz} \mapsto h_{1}\sqrt{dz}h_{2}\sqrt{dz} = h_{1}h_{2}dz.
    \end{equation}

    Now, assume that $\Omega_{1}$ and $\Omega_{2}$ are simply-connected domains, and assume for the moment that they do not contain $\infty$.
    Let $g: \Omega_{1} \rightarrow \Omega_{2}$ be a biholomorphism.
    We write $dg^{*}$ for the adjoint of the differential of $g$.
    We obtain a commutative diagram
    \begin{equation*}
        \begin{tikzcd}
            T^{*}_{1,0}\Omega_{2} \ar[r, "dg^{*}"] \ar[d] & T^{*}_{1,0}\Omega_{1} \ar[d] \\
            \Omega_{2} & \ar[l, "g"] \Omega_{1}
        \end{tikzcd}
    \end{equation*}
    This induces a map
    \begin{equation*}
        dg^{*}:\Omega^{1,0}(\Omega_{2}) \rightarrow \Omega^{1,0}(\Omega_{1}),\; \sigma \mapsto dg^{*} \circ \sigma \circ g.
    \end{equation*}
    We write $dz_{i} \in \Omega^{1,0}(\Omega_{i})$ for the canonical sections.
    We then have
    \begin{equation}\label{eq:derivative}
        dg^{*} (dz_{2}) = g' dz_{1}.
    \end{equation}

    One might now ask if $dg^{*}:T_{1,0}^{*}\Omega_{2} \rightarrow T_{1,0}^{*}\Omega_{1}$ admits a square root, i.e.~if there exists an associated map $\sqrt{dg^{*}}$ such that the following diagrams commute
    \begin{equation}\label{diag:SquareRoot}
        \begin{tikzcd}
            E|_{\Omega_{2}} \ar[r, "\sqrt{dg^{*}}"] \ar[d] & E|_{\Omega_{1}} \ar[d] & & E|_{\Omega_{2}} \otimes E|_{\Omega_{2}} \ar[d] \ar[rr, "\sqrt{dg^{*}} \otimes \sqrt{dg^{*}}" ] & & E|_{\Omega_{1}} \otimes E|_{\Omega_{1}} \ar[d] \\
            \Omega_{2} & \ar[l, "g"] \Omega_{1} & &  T^{*}_{1,0}\Omega_{2} \ar[rr, "dg^{*}"] & & T^{*}_{1,0}\Omega_{1}
        \end{tikzcd}
    \end{equation}
    It follows directly from \cref{eq:derivative,eq:pairing} that if $\sqrt{g'}$ is a square root of $g'$, then the rule
    \begin{equation}\label{eq:ActionOnHalfOrder}
        \sqrt{dg^{*}}\sqrt{dz_{2}} = \sqrt{g'}\sqrt{dz_{1}},
    \end{equation}
    defines a transformation $\sqrt{dg^{*}}$ such that Diagram \eqref{diag:SquareRoot} commutes.
    We moreover ask that the transformation $\sqrt{dg^{*}}$ be an isomorphism (of holomorphic vector bundles).
    This is the case if and only if the right-hand-side of \eqref{eq:ActionOnHalfOrder} defines a holomorphic section of $E|_{\Omega_{1}}$, which is the case if and only if $\sqrt{g'}$ is holomorphic.
    The assumption that $\Omega_{1}$ is simply-connected implies that $g': \Omega_{1} \rightarrow \C^{\times}$ admits exactly two holomorphic square roots, and thus we have the following result.
    \begin{lemma}
        For every biholomorphism $g$, there exist exactly two square roots of $dg^{*}$.
    \end{lemma}

    The pairs $(g, \sqrt{g'})$ are exactly the elements of the space $\Hdiffhat(\Omega_{1},\Omega_{2})$.
    Now, given $\sigma \in \Omega^{\frac{1}{2},0}(\Omega_{2})$, and $\hat{g} = (g, \sqrt{g'})$ we define
    \begin{equation*}
        \hat{g} \cdot \sigma = \sqrt{dg^{*}} \circ \sigma \circ g \in \Omega^{\frac{1}{2},0}(\Omega_{1}).
    \end{equation*}
    This is exactly the transformation behaviour of half-order differentials from \cref{def:HalfOrderDifferentials}.

    We refer to \cite[Sec.~2]{KW22} for further information on how to translate between half-order differentials (used in this paper) and spinors on the circle (as they appear in conformal field theory literature).

 \end{section}


\begin{thebibliography}{99}

 
\bibitem{Ahlfors} Ahlfors, L. Lectures on quasiconfonnal mappings, D. Van Nostrand, Princeton, New Jersey, 1966.
 \bibitem{Andersson} Andersson, J.-E.  On the degree of polynomial approximation in $E^p(D)$. J. Approx. Theory. {\bf 19} (1977) no 1, 61--68. 
 \bibitem{Barrett_Bolt_Mobius} Barrett, D.; and Bolt, M.  Cauchy integrals and M\"obius geometry of curves. Asian J. Math. {\bf 11} (2007), no. 1, 47--53.
 \bibitem{Barrett_Bolt_Laguerre} Barrett, D.; and Bolt, M. Laguerre arc length from distance functions. Asian J. Math. {\bf 14} (2010), no. 2, 213--233.
\bibitem{Bell_book}  Bell, S. R. The Cauchy transform, potential theory and conformal mapping. Second edition. Chapman \& Hall/CRC, Boca Raton, FL, 2016. 

\bibitem{BergmanSchiffer}  Bergman, S.; Schiffer, M.
  Kernel functions and conformal mapping. Compositio Math. {\bf 8}, (1951), 205--249.

\bibitem{Bilalov_Najafov} Bilalov, B. T.; and  Najafov, T. I. On basicity of systems of generalized Faber polynomials. Jaen J. Approx. {\bf 5} 1 (2013), 19--34. 

\bibitem{Calderon} Calder\'on, A. Cauchy integral on Lipschitz curves and related operators. Proc. Nat. Acad. Se., vol.74, 4, 1977, 1324--1327.
\bibitem{Cavus}
 \c{C}avu\c{s}, A.  Approximation by generalized Faber series in Bergman spaces on finite regions with a quasiconformal boundary.
 J. Approx. Theory {\bf 87} (1996), no. 1, 25--35.

\bibitem{Cavus_Israfilov} \c{C}avu\c{s}, A. and Israfilov, D. M. Approximation by Faber-Laurent rational functions in the mean of functions of class $L_p(\Gamma)$ with $1<p<\infty$.  Approximation Theory and its Applications, {\bf 11} (1995), no. 1, 105--118.    

\bibitem{Coifman-Fefferman}
Coifman, R and Fefferman, C. Weighted norm inequalities for maximal functions and singular integrals. Studia Math. {\bf 51} (1974), 241--250.

 \bibitem{CoifmanMeyer} Coifman, R; Meyer, Y.
Wavelets.
Calder\'on-Zygmund and multilinear operators.Cambridge Stud. Adv. Math., 48
Cambridge University Press, Cambridge, 1997. xx+315 pp.
\bibitem{David} David, G. Op\'erateurs int\'egraux singuliers sur certaines courbes du plan complexe.  Annales scientifiques de l’\'E.N.S. 4 e s\'erie, tome 17, no 1 (1984),  157--189.

\bibitem{Duren_Hp_book}  Duren, P. Theory of $H^p$ spaces.  Pure and Applied Mathematics, Vol. 38 Academic Press, New York-London 1970.

\bibitem{Duren_book} Duren, P.  Univalent functions. Grundlehren der Mathematischen Wissenschaften [Fundamental Principles of Mathematical Sciences], {\bf259}. Springer-Verlag, New York, 1983. 

\bibitem{Hawley_Schiffer_half-order}
Hawley, N. S.; and Schiffer, M.
Half-order differentials on Riemann surfaces.
Acta Math. {\bf 115} (1966), 199--236.

\bibitem{Huang}
Y.-Z. Huang, {\it Two-dimensional conformal geometry and vertex operator
  algebras}, Progress in Mathematics, vol. 148, Birkh\"auser, Boston, MA, 1997.

 \bibitem{Ibramamed} Ibragimov, I. I.; Mamedhanov, D. Constructive characterization of some classes of functions (Russian), Dokl. Akad. Nauk SSSR {\bf 223} (1975), 35--37.

\bibitem{KY2}
A.~A.~Kirillov and D.~V.~Yuriev, {\it Representations of the Virasoro
algebra by the orbit method}, J. Geom. Phys. {\bf 5} (1988), no.
3, 351--363.
 
\bibitem{Kokilasvili} Kokila\v{s}vili, V. M.  A direct theorem on mean approximation of analytic functions by polynomials, Dokl. Akad. Nauk SSSR. {\bf 185} (1969), 749--752; Soviet Math. Dokl. {\bf 10} (1969), 411--414.


\bibitem{KW22} Kristel, P.; Waldorf, K. Fusion of Implementers for spinors on the circle, Adv.~Math. {\bf{402}} (2022).

\bibitem{KS23} Kristel, P.; Schippers, E. Grassmannians of Lagrangian Polarizations, arXiv:2304.10774

\bibitem{Lehto}
Lehto, O. {Univalent functions and {T}eichm\"uller spaces}. Graduate Texts in Mathematics, Vol. 109, Springer-Verlag, New York, 1987.

\bibitem{Nagbook} Nag, S.  {The Complex Analytic Theory of
 {T}eichm\"uller Spaces}, Canadian Mathematical Society Series of
  Monographs and Advanced Texts.  John Wiley \& Sons, Inc., New York, 1988.

\bibitem{NagSullivan}  S. Nag and D. Sullivan, {\it Teichm\"uller theory and the universal period mapping via quantum calculus and the $H^{1/2}$ space on the circle}, Osaka J. Math. {\bf 32} (1995), no. 1, 1--34.

\bibitem{Nap_Yulm_Cauchy} Napalkov, V.~V., Jr.; Youlmukhametov, R. S. Criterion of surjectivity of the Cauchy transform on a Bergman space.  Israel. Math. Conf. Proc. {\bf 15}, Entire functions in modern analysis, Tel-Aviv, 1997.

\bibitem{Plymen_Robinson} Plymen, R.~J.; Robinson, P.~L., Spinors in Hilbert space. Cambride Tracts in Mathematics, Vol.~114, 1994.

\bibitem{Nap_Yulm} Napalkov, V.~V., Jr.; Yulmukhametov, R. S. On the Hilbert transform in the Bergman space. (Russian) Mat. Zametki {\bf 70} (2001), no. 1, 68--78; translation in Math. Notes {\bf 70} (2001), no. 1-2, 61--70.

\bibitem{Pommerenkebook}  Pommerenke, C. {{ Univalent functions}}, With a
  chapter on quadratic differentials by Gerd Jensen. Studia
  Mathematica/Mathematische Lehrb\"ucher, Band XXV. Vandenhoeck \& Ruprecht,
  G\"ottingen  1975.

\bibitem{FM Riesz} Riesz F.;  Riesz, M.  \"Uber die randwerte einer analtischen funktion, Compte Rendues du Qua- tri\`eme Congr\`es des Math\'ematiciens Scandinaves, Stockholm 1916, Almqvist and Wiksell, Upsala, 1920. 3

\bibitem{Schiffer_first}  Schiffer, M. The kernel function of an orthonormal system. Duke Math. J. {\bf 13} (1946). 529--540.

\bibitem{Schiffer_expository} Schiffer, M. Fredholm eigenvalues and Grunsky matrices.
   Ann. Polon. Math. {\bf 39} (1981), 149--164.

\bibitem{Schippers_Staubach_Grunsky_expository}  Schippers, E.; and Staubach, W.  Analysis on quasicircles-A unified approach through transmission and jump problems. EMS Surv. Math. Sci. {\bf 9} (2022), no.1, 31–-97.

\bibitem{WP_Thurston_Schippers_Staubach}  Schippers, E. and Staubach, W. Weil-Petersson Teichmüller theory of surfaces of infinite conformal type (2023). Accepted for publication in the collection "In the Tradition of Thurston" Volume III.

\bibitem{Segal04} Segal, G. The definition of conformal field theory. In Topology Geometry and quantum field theory, vol.~308 of London Math.~Soc.~Lecture Note Ser. (2004)



\bibitem{SegalUnitary} G. Segal, {\it Unitary representations of some infinite-dimensional groups}. Comm. Math. Phys. {\bf 80} (1981), no. 3, 301–342.


\bibitem{Shale} D. Shale, {\it Linear symmetries of free boson fields}. Trans. Amer. Math. Soc. {\bf 103} (1962), 149--167.

\bibitem{Shen} Shen, Y.
On Grunsky operator.
Sci. China Ser. A {\bf 50} (2007), no.12, 1805--1817.

\bibitem{ShenFaber} Shen, Y.  Faber polynomials with applications
to univalent functions with quasiconformal extensions,  Sci. China Ser. A {\bf 52} (2009), no. 10, 2121–-2131.

\bibitem{Shirazi_Grunsky} Shirazi, M. Faber and Grunsky Operators Corresponding to Bordered Riemann Surfaces, Conform. Geom. Dyn. {\bf 24} (2020), 177--201.

\bibitem{TNT} Takhtajan, L.; Teo, L-P. Weil-Petersson Metric on the Universal Teichm\"uller Space,
 Memoirs of the American Mathematical Society. {\bf 183}  (2006) no. 861.

 \bibitem{Tener17} Tener, J. Construction of the unitary free fermion Segal CFT. Comm.~Math.~Phys. Vol.~355 (2017)

\bibitem{Zinsmeister}
Zinsmeister, M.
Repr\'esentation conforme et courbes presque lipschitziennes.
Ann. Inst. Fourier (Grenoble) {\bf{34}}   (1984), no.2, 29–-44.

\end{thebibliography}
\end{document}